\documentstyle[10pt]{article}
\newcommand{\rom}{\mathrm}

\def\N{{\rom I\kern-.1567em N}}
\def\R{{\rom I\kern-.1567em R}}
\def\C{{\rom C\kern-6.5pt  
          \vrule height 7.7pt width 0.4pt depth -0.5pt \phantom {.}\;}}             

\newtheorem{satz}{Theorem}
\newtheorem{lemma}[satz]{Lemma}
\newtheorem{sublemma}[satz]{Sublemma}
\newtheorem{proposition}[satz]{Proposition}
\newtheorem{coro}[satz]{Corollary}

\newtheorem{remark}[satz]{Remark}
\newtheorem{defi}[satz]{Definition}
\newcommand{\eps}{\varepsilon}
\newcommand{\falle}{\;\;\forall}
\newcommand{\gen}{\rightarrow}
\newcommand{\nach}[1]{\stackrel{#1}\rightarrow}
\newcommand{\Norm}[1]{\Bigl\|#1\Bigr\|}
\newcommand{\norm}[1]{\|#1\|}
\newcommand{\Betr}[1]{\Bigl| #1  \Bigr|}
\newcommand{\betr}[1]{| #1  |}
\newcommand{\eing}[1]{|_{#1}}
\newcommand{\Grabstein}{\rule{1.2ex}{1.2ex}}
\newcommand{\ebew}{\hfill\Grabstein}
\newcommand{\bgl}{\begin{eqnarray}}
\newcommand{\bglst}{\begin{eqnarray*}}
\newcommand{\egl}{\end{eqnarray}}
\newcommand{\eglst}{\end{eqnarray*}}
\newcommand{\Pel}{Pe\l\-czy\'ns\-ki}
\newcommand{\folgt}{\Rightarrow}
\newcommand{\leins}{l^{1}}
\newcommand{\lunendl}{l^{\infty}}
\newcommand{\Ref}[1]{(\ref{#1})}
\newcommand{\gdw}{\Leftrightarrow}
\newcommand{\eins}{{1}}
\newcommand{\id}{{{\rom i}\rom{d}}}
\newcommand{\Lst}{$Li^*$}
\newcommand{\lst}{{li}^*}
\newcommand{\genwst}{\nach{w^{*}}}

\newcommand{\gentaumu}{\nach{\tau_\mu}}
\newcommand{\genw}{\nach{w}}
\newcommand{\tauseq}{\tau_{\lst}}
\newcommand{\gentauseq}{\nach{\tauseq}}
\newcommand{\wst}{$w^{*}$}
\newcommand{\w}{$w$}
\newcommand{\ball}[1]{{\rom B}_{#1}}
\newcommand{\mdE}{|\;}
\newcommand{\Lunendl}{{{\rom L}^{\infty}}}
\newcommand{\Leins}{{\rom L}^{1}}

\newcommand{\LLa}{\raisebox{-1.1ex}{'\kern-0.09em'}}
\newcommand{\Ree}{{\rom R}{\rom e}\,}
\newcommand{\asyleins}{\stackrel{{\rom{a}\rom{s}\rom{y}}}{\sim}\leins}
\newcommand{\almleins}{\stackrel{{\rom{a}\rom{l}\rom{m}}}{\sim}\leins}
\newcommand{\isoleins}[1]{\stackrel{#1}{\sim}\leins}
\newcommand{\dist}{\mbox{dist}}
\begin{document}
\begin{center}
{\bf L-embedded Banach spaces and measure topology}
\medskip\\ H. Pfitzner\end{center}
\begin{abstract}\noindent
An L-embedded Banach spaace is a Banach space which is complemented in its
bidual such that the norm is additive between the two complementary parts.
On such spaces we define a topology, called an abstract measure
topology, which by known results coincides with the usual measure
topology on preduals of finite von Neumann algebras (like $\Leins([0,1])$).
Though not numerous, the known properties of this topology suffice
to generalize several results on subspaces of
$\Leins([0,1])$ to subspaces of arbitrary L-embedded spaces.
\end{abstract}
\bigskip\bigskip\bigskip
\noindent
{\bf\S 1 Introduction}\\
This article continues the investigations made in \cite{Pfi-Fixp, Pfi-vN}
on asymptotically isometric copies of $\leins$ in preduals of von Neumann
algebras and in L-embedded Banach spaces.
(For defintions see below.)
In \cite{Pfi-vN} it has been proved that, roughly speaking, in the predual
of a finite von Neumann algebra the only non-trivial bounded
sequences that converge to $0$ with respect to the measure topology
are essentially those that span $\leins$ asymptotically;
for $\Leins(\mu)$, $\mu$ a finite measure, this characterization
has been known for quite a time
\cite[Th.\ 2]{KadPel}, \cite[Th.\ 3, Rem.\ 6bis]{Pi-bases}.

From the point of view of Banach space theory, L-embedded Banach spaces
provide a natural frame for preduals of von Neumann algebras.
So the starting point of this paper is on the one hand
the definition of an abstract measure topology, Definition \ref{defMass},
patterned after the just mentionend characterization and on the other hand
the easy but important observation, Theorem \ref{theo_existence},
that every L-embedded space admits such a topology.
Although this topology does not come out easily with its properties 
- at the time of this writing it is not clear whether it is
Hausdorff let alone metrizable or whether addition is continuous -
it allows to generalize several results on subspaces of
$\Leins(\mu)$ to subspaces of arbitrary L-embedded spaces.
Thus section \S 4 of the present paper is titled "Section IV.3 of
\cite{HWW} (partly) revisited".
For example, Theorem \ref{theoBL} generalizes a theorem of
Buhvalov-Lozanovskii which
describes the link between L-embeddedness and measure topology
for subspaces $Y$ of $\Leins(\mu)$, $\mu$ finite:
$Y$ is L-embedded if and only if its unit ball is closed in measure.
(Note in passing that this criterion involves only the space $Y$
itself, not its bidual.) Moreover, as a consequence of this,
the closedness in measure of the unit ball of $Y$ is a weak
substitute for compactness which could be called "convex sequential
compactness", see Corollary \ref{coro_compakt}.
We also reprove a result of Godefroy and Li
concerning a criterion for L-embedded subspaces
which are duals of M-embedded spaces, see Theorem \ref{theo3.10}.
In this vein, that is by substituting arbitrary L-embedded spaces for
$\Leins(\mu)$, we recover also some results of Godefroy, Kalton, Li
\cite{GoKaLi} in \S 5.
Finally, in \S 6 it is proved that addition is $\tau_{\mu}$-continuous
in preduals of von Neumann algebras.
\bigskip\bigskip\\
{\bf\S 2 Notation, Background}:\\
The results are stated for complex scalars.
The dual of a Banach space $X$
is denoted by $X'$.
$\ball{X}$ denotes the unit ball of $X$.
Subspace of a Banach space means norm-closed subspace,
bounded always means norm-bounded.
As usual, we consider a Banach
space as a subspace of its bidual omitting the canonical embedding.
$[x_n]$ denotes the closed linear span of a (finite or infinite) sequence
$(x_n)$.

Basic properties and definitions which are not explained here
can be found in \cite{Die-Seq} or in \cite{LiTz1}-\cite{LiTz2}
for Banach spaces and in \cite{Ped, Tak} for C$^*$-algebras.
The standard reference for M- and L-embedded spaces is the monograph
\cite{HWW}.\medskip\\
Let $(x_n)$ be a sequence of nonzero elements in a Banach space $X$.

We say that
\begin{em}$(x_n)$ spans $\leins$ isomorphically\end{em}
(or \begin{em}$r$-isomorphically\end{em} to be more precise)
- $(x_n)_{n\in\N}\isoleins{r}$ or just
$x_n\isoleins{r}$ in symbols -
if there exists $r>0$ (trivially $r\leq 1$) such that
$r(\sum_{n=1}^{\infty}\betr{\alpha_n})\leq
\norm{\sum_{n=1}^{\infty} \alpha_n \frac{x_n}{\norm{x_n}}}   \leq
        \sum_{n=1}^{\infty}\betr{\alpha_n}$
for all scalars $\alpha_n$ (the second inequality being trivial).

We say that
\begin{em}$(x_n)$ spans  $\leins$ almost isometrically\end{em}
- $x_n\almleins$ in symbols -
if there is a sequence $(\delta_m)$ in $[0,1[$
tending to $0$ such that
$(x_n)_{n\geq m}\isoleins{1-\delta_m}$ for all $m\in\N$.
Recall that the Banach-Mazur distance of two Banach spaces $X$ and $Y$
is defined by $\dist(X,Y)=\inf\norm{T}\,\norm{T^{-1}}$
where the infimum extends over all
surjective isomorphisms $T:X\gen Y$.
To avoid confusion, notice that $\dist(\leins, [x_n])=1$
is not the same as $x_n\almleins$.

Finally, a sequence $(x_n)$ is said {\em to span $\leins$
asymptotically isometrically} or just
{\em to span $\leins$ asymptotically}
- $x_n\asyleins$ in symbols -
if there is a sequence
$(\delta_n)$ in $[0,1[$ tending to $0$ such that
$\sum_{n=1}^{\infty}(1-\delta_n)\betr{\alpha_n}\leq
\norm{\sum_{n=1}^{\infty} \alpha_n \frac{x_n}{\norm{x_n}}}   \leq
         \sum_{n=1}^{\infty}\betr{\alpha_n}$
for all scalars $\alpha_n$.

Note that the present definitions of almost and asymptotically isometric
differ slightly from those in \cite{DJLT}, \cite{Pfi-Fixp} by the term
$x_n /\norm{x_n}$ but that, of course, for normalized sequences
the definitions are the same.
Note also the technical detail that because of this term
one might have $\norm{x_n}\gen 0$ for a sequence spanning $\leins$
isomorphically
whereas sequences that are
equivalent to the canonical $\leins$-basis (\cite[p.\ 43]{Die-Seq})
are uniformly bounded away from $0$.
We say that a Banach space is isomorphic
(respectively almost isometric respectively asymptotically isometric)
to $\leins$ if it has a basis with
the corresponding property.\medskip\\
Let $Y$ be a subspace of a Banach space $X$ and $P$ be a projection
on
$X$.
$P$ is called an \begin{em}L-projection\end{em}
provided
$\norm{x}=\norm{Px}+\norm{(\id_{X}-P)x}$ for all $x\in X$.
A subspace $Y\subset X$ is called an \begin{em}M-ideal in $X$\end{em}
if its annihilator $Y^{\bot}$ in $X'$ is the
range of an L-projection on $X'$.
$Y$ is called an \begin{em}L-summand in $X$\end{em} if it is the
range of an
L-projection on $X$.
In the special case in which $X=Y''$ and in which $Y$ is an M-summand
(respectively an L-summand) in $Y''$  we say
that $Y$ is \begin{em}M-embedded\end{em} (respectively
\begin{em}L-embedded\end{em}).
As examples we only mention that preduals of von Neumann algebras,
in particular $\leins$ and $L^1$-spaces,
furthermore the Hardy space $H_0^1$ and the dual
of the disc algebra are L-embedded.
The sequence space $c_0$,
the space of compact operators on a Hilbert space, and
the quotient $C/A$ of the
continuous functions on the unit circle by the disc algebra $A$
are examples among M-embedded spaces.
It is not difficult but important to see that if there is an L-projection $P$
on a Banach space $X$ then each contractive projection on $X$
which has the same kernel as $P$ coincides with $P$,
see \cite[Prop.\ I.1.2]{HWW}.
It follows that if $X$ is M-embedded then the canonical inclusion of
$X'$ in $X'''$ is an L-summand in $X'''$ that is $X'$ is L-embedded;
the converse is false \cite[III.1.3]{HWW};
in fact, for an L-embedded space being the dual of an M-embedded space
can be quite a restrictive condition: For instance, while the dual of
any C$^*$-algebra is L-embedded only those C$^*$-algebras are M-embedded
which are isometrically $^*$-isomorphic to the algebra of compact operators
or to a $c_0$-sum of such algebras \cite[III.2.9]{HWW}.
Throughout this note, if $X$ denotes an L-embedded Banach space (which is not
always the case) we will write $X_s$ for
the complement of (the canonical embedding of) $X$ in $X''$
that is $X''=X\oplus_1 X_s$.
In this case $P$ will denote the L-projection from $X''$ onto $X$.\\
We recall Godefroy's fundamental result \cite{G-bien}, \cite[IV.2.2]{HWW}
that L-embedded Banach spaces are $w$-sequentially complete.
This will be used mostly without reference; together with Rosenthal's
$\leins$-theorem a typical application is that each bounded sequence in
an L-embedded space contains a subsequence which either spans $\leins$
or converges weakly.
There is a useful criterion for L-embeddedness of
subsapces of L-embedded spaces due to
Li (\cite{Li-Ox} or \cite[Th.\ IV.1.2]{HWW}) which we state for easy
reference:
\begin{lemma}[Li]\label{lem_Li}
For an L-embedded Banach space $X$
(with L-decomposition $X''=X\oplus_1 X_s$ and L-projection $P$
on $X''$ with range $X$)
and a closed subspace $Y$ of $X$ the following assertions are equivalent.
\begin{description}
\item{\makebox[1.8em]{(i)}}
$Y$ is L-embedded.
\item{\makebox[1.8em]{(ii)}}
$Y^{\bot\bot}=Y\oplus_1 (Y^{\bot\bot}\cap X_s)$.
\item{\makebox[1.8em]{(iii)}}
$P\overline{\ball{Y}}^{w^*}=\ball{Y}$.
\item{\makebox[1.8em]{(iv)}}
$P Y^{\bot\bot}=Y$.
\end{description}
In particular if $Y$ is L-embedded and if one identifies
$Y''=Y\oplus_1 Y_s$ and $Y^{\bot\bot}\subset X''$
then $Y_s=Y^{\bot\bot}\cap X_s$.
\end{lemma}
Let us finally cite some technical results from \cite{Pfi-Fixp} which will be used
in the sequel.\\
It is routine to show that sequences spanning $\leins$
asymptotically are stable by adding norm-null sequences
\cite[Lem.\ 4]{Pfi-vN}, to be more precise,
\begin{em}let $(x_n)$, $(y_n)$ be two sequences in a Banach space $X$ such
that $(x_n)$ spans $\leins$ asymptotically, $\inf\norm{x_n}>0$,
$\norm{y_n}\gen 0$ and $x_n +y_n\neq 0$.
Then $(x_n + y_n)$ spans $\leins$ asymptotically, too.\end{em}\\
Although it has been proved in \cite{DJLT} that there are 
almost isometric $\leins$-copies which do not contain asymptotic
ones, both notions "coincide up to subsequences" in L-embedded Banach spaces,
more precisely, in L-embedded Banach spaces each sequence
spanning $\leins$ almost isometrically contains a subsequence
spanning $\leins$ asymptotically \cite[Cor.\ 3]{Pfi-Fixp}.\\
The following lemma is fundamental for the rest of the paper.
It is  an immediate consequence
of \cite{Pfi-Fixp}
and says that within L-embedded spaces, sequences spanning
$\leins$ almost or asymptotically isometrically behave like the standard
basis of $\leins$ as to their \wst-accumulation points.
Mostly it will be used with $M$ being a countable set of normalized
elements that span $\leins$ asymptotically.
\begin{lemma}\label{lem_HP}
Let $X$ be L-embedded (with L-decomposition $X''=X\oplus_1 X_s$),
let $M\subset X$ be a subset of the unit sphere of $X$.
Then the following assertions are equivalent.
\begin{description}
\item{\makebox[1.8em]{(i)}}
$M$ contains a sequence spanning $\leins$ asymptotically.
\item{\makebox[1.8em]{(ii)}}
$M$ admits a \wst-accumulation point $x_s\in X_s$ of norm one.
\end{description}
In fact, each \wst-accumulation point of a sequence spanning $\leins$
asymptotically lies in $X_s$.
\end{lemma}
{\em Proof}:
(i)$\folgt$(ii) and the second statement:
see last paragraph of the proof of
\cite[Lem.\ 1]{Pfi-Fixp};
(ii)$\folgt$(i): \cite[Th.\ 2]{Pfi-Fixp}\ebew
\bigskip\bigskip\\
{\bf\S 3 {Abstract measure topology}}\\
To prepare the definiton of an abstract measure topology we recall some facts
about sequential spaces (see for example
\cite[1.6-1.7]{Engel2} or \cite{Kisynski}, cf. also \cite{Kelley})
because the topology will be defined by
determining the class of its convergent sequences.

A topological space is called a sequential
space if closedness and sequential closedness coincide.
A topological space is called a Fr\'echet space if closure and
sequential closure coincide.
Clearly first countable spaces are Fr\'echet spaces
and Fr\'echet spaces are sequential spaces.
Recall that a topological space is a $\mbox{T}_1$-space if every one-point set
is closed; this happens if each convergent sequence has a unique limit.

A \Lst-space\footnote{To avoid confusion with the letter L like in
L-embedded, L-projection, L-structure we prefer the notation
\Lst instead of ${\cal L}^*$ as in \cite{Engel2, Kisynski}.}
is a triple $(X, {\cal C}, \lst)$
where $X$ is a set,
${\cal C}\subset X^{\N}$
a class of sequences of $X$ (called the convergence class)
and $\lst:{\cal C}\rightarrow X$
a map (called limit operator) satisfying the following conditions
(L1)-(L3).
We write $\lst x_n$ instead of $\lst((x_n)_{n\in\N})$;
the elements of ${\cal C}$ are called ${\cal C}$-convergent
sequences.\\
Let $(x_n)$ be any sequence in $X$, let $x\in X$.
\begin{description}
\item{\makebox[1.3em]{(L1)}}$\;\;$ If $x_n=x$ for all $n\in\N$ then $(x_n)\in{\cal C}$
       and $\lst x_n=x$.
\item{\makebox[1.3em]{(L2)}}$\;\;$ If $(x_n)\in{\cal C}$ with
    $\lst x_n=x$ then $(x_{n_k})\in{\cal C}$ and
     $\lst x_{n_k}=x$ for each subsequence $(x_{n_k})$ of $(x_n)$.
\item{\makebox[1.3em]{(L3)}}$\;\;$ If a sequence $(x_n)$ is such that there is $x\in X$
    and any subsequence $(x_{n_k})$
   contains a further subsequence $(x_{n_{k_m}})$ such that
       $(x_{n_{k_m}})\in{\cal C}$ and $\lst x_{n_{k_m}}=x$
   then $(x_n)\in{\cal C}$ and $\lst x_n=x$.
\end{description}
On a \Lst-space $(X, {\cal C}, \lst)$ one defines
a topology $\tauseq$ - called the sequential topology
induced by $\lst$ - by taking as the family of closed sets all
$\lst$-sequentially closed sets; here we call a set $A$
$\lst$-sequentially closed if
$\lst x_n \in A$ for all ${\cal C}$-convergent sequences $(x_n)$
contained in $A$.
It is elementary to verify that in this way one indeed
obtains a topology and that the $\tauseq$-convergent sequences
are exactly the ${\cal C}$-convergent sequences,
that is for every sequence $(x_n)$ in $X$
one has $\lst x_n =x$ if and only if $x_n \gentauseq x$.\smallskip

An ${\cal S}^*$-space is a \Lst-space $(X, {\cal C}, \lst)$ satisfying
additionally
\begin{description}
\item{\makebox[1.3em]{(L4)}}$\;\;$
If $(x_m^{(n)})_{m}\in{\cal C}$ for all $n\in\N$
and $(x_n)\in{\cal C}$ such that $\lst x_m^{(n)}=x_n$ for all $n\in\N$
and $\lst x_n=x$ then there exist two sequences $(n_k)$, $(m_k)$ such
that $\lst x_{m_k}^{m_k}=x$.\smallskip
\end{description}
Endowed with $\tauseq$, an ${\cal S}^*$-space becomes a Fr\'echet space
\cite[1.7.18,19]{Engel2}.\smallskip

As already mentioned in the introduction the following definition is
patterned after a characterization of bounded measure-null sequences
in the preduals of finite von Neumann algebras \cite[Th.\ 1]{Pfi-vN}.
\begin{defi}\label{defMass}
Let $X$ be a Banach space.
A system $\tau_{\mu}$ of subsets  of $X$ is called an abstract measure topology
if it satisfies the following four conditions.
\begin{enumerate}
\item $(X,\tau_{\mu})$ is a sequential space in which every convergent
       sequence has a unique limit.
\item $\tau_{\mu}$ is weaker than the norm topology.
\item $\tau_{\mu}$ is translation invariant for sequences more precisely,
$x_n\gentaumu x$ if and only if $x_n -x\gentaumu0$ for any sequence
$(x_n)$ in $X$.
\item Each bounded
sequence in $X$ that spans $\leins$ asymptotically
$\tau_{\mu}$-converges to $0$,\\
and each sequence in $X$ that $\tau_{\mu}$-converges to $0$ is bounded
and contains a subsequence which spans $\leins$ asymptotically
or tends to $0$ in norm.
\end{enumerate}
\end{defi}\bigskip

The following result is quite easy to prove. Nevertheless,
because of its importance, we call it a theorem.
\begin{satz}\label{theo_existence}
Every L-embedded Banach space admits an abstract measure topology.
\end{satz}
{\em Proof}:
Let $X$ be an L-embedded Banach space with L-decomposition
$X''=X\oplus_1 X_s$. Set
\bgl
{\cal C}_0 =\{(x_n)\mdE&&\!\!\!\!\!\!\!\!\!\!\! (x_n) \mbox{ is bounded and
                 every subsequence }(x_{n_k})\mbox{ of }(x_n)
               \mbox{ contains}
\nonumber\\
&& \mbox{ a subsequence }(x_{n_{k_l}})
             \mbox{ such that } x_{n_{k_l}}\asyleins \mbox{ or }
               \norm{x_{n_{k_l}}}\gen 0\},
\nonumber\\
{\cal C} =\{(x_n)\mdE&&\!\!\!\!\!\!\!\!\!\!\! \mbox{ there exists }x\in X \mbox{ such that }
                 (x_n -x)\in {\cal C}_0\}.
\nonumber
\egl
We define a limit operator $\lst:{\cal C} \gen X$ by $\lst x_n=x$ where $x\in X$
is such that $(x_n -x)\in {\cal C}_0$.
To show that $(X, {\cal C}, \lst)$ is a \Lst-space the only thing to verify
is that $\lst$ is well defined as a map.
because then conditions (L1) - (L3) are immediate
from the definiton of ${\cal C}$.

Suppose that there are $x,y\in X$, $(x_n)\in {\cal C}$ such that both
$(x_n -x)\in {\cal C}_0$ and $(x_n -y)\in {\cal C}_0$. If $(x_n -x)$ or $(x_n -y)$
admits a subsequence tending to $0$ in norm then $x=y$ because the
norm topology is Hausdorff.
Otherwise, after passing to an appropriate subsequence,
we suppose that both sequences are uniformely bounded away from $0$ in norm
and that both $x_n-x \asyleins$ and $x_n-y \asyleins$.
Since both sequences are bounded,
by Lemma \ref{lem_HP}
they admit two \wst-accumulation points $x_s, y_s \in X_s$
and there is a net $(x_{n_\gamma})$ such that
$x_{n_\gamma} -x\genwst x_s$ and
$x_{n_\gamma} -y\genwst y_s$. But this means that
$x_{n_\gamma}\genwst x+x_s =y+y_s$ whence $x=y$ (and $x_s =y_s$).

We define the abstract measure topology $\tau_{\mu}$ as the
sequential topology induced by $\lst$. It is immediate from the definition
of ${\cal C}$ that $\tau_{\mu}$ satisfies the conditions of
Definition \ref{defMass}.
\ebew\bigskip\\
Lemma \ref{lem_elementar} shows some elementary properties of $\tau_{\mu}$
which in the sequel will be used mostly without reference.

\begin{lemma}\label{lem_elementar}
If a Banach space $X$ admits an abstract measure
topology $\tau_{\mu}$ then $\tau_{\mu}$ has the following properties.
\begin{description}
\item{(a)} $(X,\tau_{\mu})$ is a $\mbox{T}_1$-space.
\item{(b)}
The relative topology of $\tau_{\mu}$ on a subspace of $X$ is again
an abstract measure topology.
\item{(c)} Closedness and sequential closedness coincide for $\tau_{\mu}$.
Sequentially continuous maps on $X$ are continuous.
\item{(d)}
If $X$ does not contain a copy of $\leins$ the norm topology is an abstract
measure topology and is the only one.
\item{(e)} $\tau_{\mu}$ is unique.
\item{(f)}
If $X$ is the predual of a finite von Neumann
algebra then $\tau_{\mu}$ coincides on bounded sets with the usual
measure topology; on unbounded sets it does not in general.
\item{(g)} Multiplication by scalars is $\tau_{\mu}$-continuous.
\item{(h)} If $X$ is L-embedded (i.e. $X''=X\oplus_1 X_s$)
and if a net $(x_{\gamma})$ in $X$ 
\wst-converges to $x_s\in X_s$ such that
$\norm{x_{\gamma}}\leq\norm{x_s}$ then $x_{\gamma}\gentaumu 0$.
\item{(i)} If $X$ is L-embedded then for any $x''\in \overline{\ball{X_s}}^{w^*}$
there exists a net $(x_{\gamma})$ in $\ball{X}$ such that both
$x_{\gamma}\gentaumu 0 \mbox{ and } x_{\gamma}\genwst x''$.
\end{description}
\end{lemma}
{\em Proof}:
(a) - (d) are clear from the definition.\\
(e) The topology of sequential spaces is determined by its
convergent sequences. Thus the abstract measure topology
is unique because the conditions of Definition \ref{defMass}
determine all convergent sequences.\\
(f) The assertion concerning bounded sets
follows from \cite[Th.\ 1]{Pfi-vN}.
For the assertion concerning unbounded sets we consider
the usual measure (=pointwise) topology and the \wst-topology on $\leins$:
These two and $\tau_{\mu}$ coincide on bounded sets.
But the unbounded sequence $(n\, e_n)$ converges in the usual measure
topology while it does not with respect to $\tau_{\mu}$.
(Here $(e_n)$ denotes the standard basis of $\leins$.)\\
(g) Let $\lambda_n\gen\lambda$ in $\C$,  and $x_n\gentaumu x$ in $X$.
If $\lambda_n\gen 0$ or $\norm{x_n-x}\gen0$ then
$\lambda_n x_n\gen\lambda x$ with respect
to the norm-topology and thus also with respect to $\tau_{\mu}$.
If $\inf\betr{\lambda_n}> 0$, $x_n -x\asyleins$,
and $\inf\norm{x_n -x}>0$  then
$\lambda_n x_n -\lambda_n x\asyleins$ whence
$\lambda_n x_n -\lambda x\asyleins$ because
$\norm{\lambda_n x-\lambda x}\gen0$ and because sequences spanning $\leins$
asymtotically are stable by adding norm-null sequences.
Hence $\lambda_n x_n \gentaumu \lambda x$.
Up to the usual reasoning on subsequences of subsequences this proves
that
$\lambda_n \gen \lambda$ and $x_n \gentaumu x$
imply $\lambda_n x_n\gentaumu \lambda x$.\\
(h) We assume to the contrary that $(x_{\gamma})$ does not
$\tau_{\mu}$-converge to $0$. Then there is a $\tau_{\mu}$-neighborhood 
${\cal O}$ of $0$
and a subnet $(x_{\gamma'})$ which does not
meet ${\cal O}$ and which still \wst-converges to $x_s$.
But by Lemma \ref{lem_HP} there is a sequence
$(x_{\gamma'_n})$ that spans $\leins$ almost isometrically i.e.
$x_{\gamma'_n}\gentaumu0$ hence $(x_{\gamma'_n})$ does meet ${\cal O}$.\\
(i) Let ${\cal U}$ be a \wst-neighbourhood of $x''$. Then there exists a
$x_s \in \ball{X_s} \cap {\cal U}$ and ${\cal U}$ is also a \wst-neighbourhood of $x_s$.
Let $(x_{\gamma})\subset \ball{X}$ be a net with \wst-limit $x_s$ and
such that $\norm{x_{\gamma}}=\norm{x_s}$. Let ${\cal V}$ be a $\tau_{\mu}$-neighbourhood
of $0$. Then $x_{\gamma}\gentaumu 0$ by (h) and there is $\gamma_0$ such that
$x_{\gamma}\in {\cal U}$  and $x_{\gamma}\in {\cal V}$ for all $\gamma\succeq\gamma_0$.
That is ${\cal V} \cap {\cal U} \cap \ball{X}\neq \emptyset$ whence the assertion.
\ebew\bigskip\\
Remarks:\\
1. Part (h) of the Lemma above corresponds to \cite[IV.3.7]{HWW},
part (i) to \cite[Lem.\ 2.2]{GoKaLi}.
Since in general $\ball{X_s}$ is not \wst-closed there are $\tau_{\mu}$-null sequences
which differ from the ones described in part (h). It is tempting to suppose
that the $\tau_{\mu}$-null nets are exactly those that
admit \wst-limits in $\overline{\ball{X_s}}^{w^*}$
but at the time of this writing this is not at all clear.\\
2. There exists a non-Hausdorff Fr\'echet space in which
every sequence has at most one limit \cite[1.6E]{Engel2}.
Therefore the question whether $\tau_{\mu}$ is Hausdorff is not
necessarily trivial.\\
3. If $X, Y$ are two L-embedded Banach spaces, $Y\subset X$ a subspace of $X$,
then by (b) and (e) one can identify the intrinsic abstract measure topology
of $Y$ with the relative topology of the abstract measure topology of $X$.
Therefore, in theorems like Theorem \ref{theo3.10} or Proposition
\ref{theo2.1a} it is no longer necessary to consider a sourrounding
L-embedded Banach space like $\Leins(\mu)$ in the corresponding theorems
\cite[IV.3.10]{HWW} or \cite[Prop.\ 2.1]{GoKaLi}.
This observation sheds also some light on the remarks after \cite[IV.3.5]{HWW}
and after \cite[Def.\ IV.4.2]{HWW} on the use of "nicely placed"
and "L-embedded".\\
4. It might be usefull for other purposes to modify Definition \ref{defMass}.
For example one could replace "asymptotically" by "almost isometrically"
in Definition \ref{defMass}; by \cite[Cor.\ 3]{Pfi-Fixp} this would give the
same topology for L-embedded spaces.
As a less trivial modification one could first restrict Definition \ref{defMass} to
bounded subsets of a Banach space $X$ and then define the abstract measure topology
on the whole of $X$ as the inductive limit of the
family $\tau_{\mu}\eing{n\ball{x}}$. In this case, the results of this paper would remain
valid up to some minor modifications because in all proofs except for Lemma
\ref{lem_elementar} only the restriciton of $\tau_{\mu}$ to
the unit ball is considered.
In passing we note that we restrict our attention to bounded sets mainly
because the characterization of measure-null sequences in \cite{Pfi-vN} does
not work for unbounded sequences, see the second remark after the proof of
Theorem 1 in \cite{Pfi-vN}.
\bigskip\\
We end this section with a modest attempt to get closer to the sequential
structure of $\tau_{\mu}$. In case the addition is $\tau_{\mu}$-continuous,
at least on bounded sets $\tau_{\mu}$ gives a Fr\'echet space. This applies
to von Neumann preduals which in \S 6 below will be shown to have
$\tau_{\mu}$-continuous addition.
\begin{lemma}\label{lem_Frechet}
Let $X$ be an L-embedded Banach space.
If the addition is $\tau_{\mu}$-continuous then the restriction of 
$\tau_{\mu}$ to a bounded subset of $X$ makes this set a Fr\'echet space.
\end{lemma}
{\em Proof}:
To avoid trivialities we consider a bounded sequence $x_n\asyleins$ in $X$,
and a uniformely bounded sequence of sequences $(x_m^{(n)})_{m\in\N}$ 
such that $x_m^{(n)}-x_n\asyleins$ for all $n\in\N$.
It is enough to show the existence of a sequence $(m_n)$ such that
$x_{m_n}^{(n)}\gentaumu 0$ because this will prove that $(X, {\cal C}, \lst)$
of the proof of Theorem \ref{theo_existence} satisfies (L4).

We set $y_m^{(n)}=x_m^{(n)} -x_n$. Since all $x_m^{(n)}$ are uniformely bounded there is no
loss of generality if we suppose that $\norm{y_m^{(n)}}=1$ for all $m$ and
all $n$.
By Lemma \ref{lem_HP} each universal net $(y_{m_\gamma}^{(n)})_{\gamma\in\Gamma}$
\wst-converges to a limit $y_s^{(n)}\in X_s$ of norm one.
What follows is a straightforward modification of the proof of
\cite[Th.\ 2]{Pfi-Fixp}.
Let $(\delta_n)$ be a sequence of strictly positive numbers in $]0,1]$
converging to $0$.
Set $\eta_1=\frac{1}{3}\delta_1$ and
$\eta_{n+1}=\frac{1}{3}\min(\eta_n,\delta_{n+1})$ for
$n\in\N$.
By induction over $n\in\N$ one constructs $m_n\in\N$
such that 
\bgl
\Bigl(\sum_{i=1}^{n}(1-\delta_i)\betr{\alpha_i}\Bigr)
+\eta_n \sum_{i=1}^{n}\betr{\alpha_i}
\leq
\norm{\sum_{i=1}^{n}\alpha_i y_{m_i}^{(i)}}
\mbox{ for all } n\in\N, \,\,\alpha_i\in\C.
\label{glFrechet1}
\egl
The first induction step is settled by setting $m_1=1$.
For the induction step $n\mapsto n+1$
fix an element $\alpha=(\alpha_i)_{i=1}^{n+1}$ in the
unit sphere of $\leins_{n+1}$ such that $\alpha_{n+1}\neq 0$.
The \wst-convergence (along $\gamma$) of
$(\sum_{i=1}^n \alpha_i y_{m_i}^{(i)}) + \alpha_{n+1}y_{m_\gamma}^{(n+1)}$ to
$(\sum_{i=1}^n \alpha_i y_{m_i}^{(i)}) + \alpha_{n+1}y_s^{(n+1)}$ yields
\bglst
\liminf_{\gamma}
\Norm{\Bigl(\sum_{i=1}^n \alpha_i y_{m_i}^{(i)}\Bigr) 
+ \alpha_{n+1}y_{m_\gamma}^{(n+1)}}
&\geq&
\Norm{\Bigl(\sum_{i=1}^n \alpha_i y_{m_i}^{(i)}\Bigr) 
+ \alpha_{n+1}y_s^{(n+1)}}\\
&=&\Norm{\sum_{i=1}^n \alpha_i y_{m_i}^{(i)}} + \betr{\alpha_{n+1}}\\
&\stackrel{\Ref{glFrechet1}}{\geq}&
 \Bigl(\sum_{i=1}^{n+1}(1-\delta_i)\betr{\alpha_i}\Bigr)
        +\min(\eta_n,\delta_{n+1})
\eglst
whence
\bglst
\Norm{\Bigl(\sum_{i=1}^n \alpha_i y_{m_i}^{(i)}\Bigr) 
+ \alpha_{n+1}y_m^{(n+1)}}
\geq
\Bigl(\sum_{i=1}^{n+1}(1-\delta_i)\betr{\alpha_i}\Bigr) + 3\eta_{n+1}
\eglst
for infinitely many $m$.
This extends to a finite $\eta_{n+1}$-net
$(\alpha^{l})_{l=1}^{L_{n+1}}$ in the unit sphere
of $\leins_{n+1}$ with $2\eta_{n+1}$ instead of $3\eta_{n+1}$,
to all $\alpha$ in the unit sphere of $\leins_{n+1}$
with $\eta_{n+1}$ instead of $2\eta_{n+1}$,
and finally to all $\alpha\in\leins_{n+1}$.
The details are the same as in the proof of
\cite[Th.\ 2]{Pfi-Fixp}.\\
This ends the induction. By \Ref{glFrechet1} we have
$y_{m_n}^{(n)}\gentaumu 0$. Hence $x_{m_n}^{(n)}=y_{m_n}^{(n)}+x_n\gentaumu0$ because
addition is supposed to be $\tau_{\mu}$-continuous.
\ebew\bigskip\bigskip\\
{\bf\S 4 Section IV.3 of \cite{HWW} (partly) revisited}\\
Proposition \ref{prop_reflex} generalizes some well known facts, see e.g.
\cite[Th.\ IV.8.12]{DS} for (a) and \cite[p.\ 202]{HWW} for (b).
\begin{proposition}\label{prop_reflex}
Let $X$ be an L-embedded Banach space
with its abstract measure topology
$\tau_{\mu}$.
Then the following statements hold.
\begin{description}
\item{(a)} A sequence converges in norm if (and only if) it converges
both weakly and with respect to $\tau_{\mu}$, and all limits coincide.
\item{(b)} A norm closed subspace $Y\subset X$ is reflexive if and only if
$\tau_{\mu}$ and the norm topology coincide on the unit ball of $Y$.
\end{description}
\end{proposition}
{\em Proof}:
(a) The statement is almost immediate from the definiton of $\tau_{\mu}$:
First remark that a sequence which $\tau_{\mu}$-converges to $0$
contains a subsequence which either converges to $0$ in 
norm or is uniformely bounded away from 0 in norm
and spans $\leins$ (asymptotically);
but the latter case is excluded if the sequence
also converges weakly (to whatever limit)
because $\leins$-bases do not converge weakly.
Now let $(x_n)$ be a sequence in $X$, let $x,y\in X$ be such that
both $x_n\gentaumu x$ and $x_n\genw y$. Then by what
has just been remarked, for each subsequence $(x_{n_k})$
there is a subsequence $(x_{n_{k_l}})$ such that
$x_{n_{k_l}}-x\gen 0$ in norm
for any subsequence $(x_{n_k})$ whence $x=y$ and the assertion
follows.\\
(b) By $\tau_{\norm{\cdot}}$ we denote the norm topology on $X$.
Let $Y$ be reflexive.
To show that $\tau_{\mu}$ and $\tau_{\norm{\cdot}}$ coincide on
the unit ball $\ball{Y}$ of $Y$
it is enough to show that each subsequence of a $\tau_{\mu}$-convergent
sequence in $\ball{Y}$ admits
a subsequence which converges in norm to the same limit.
But if $Y$ is reflexive then each bounded sequence contains
a weakly convergent subsequence which if the sequence is also
$\tau_{\mu}$-convergent converges in norm by (a).
Thus $\tau_{\mu}$ and $\tau_{\norm{\cdot}}$ coincide on $\ball{Y}$.\\
Conversely suppose that
$\tau_{\mu}$ and $\tau_{\norm{\cdot}}$ coincide on $\ball{Y}$.
In order to prove that $Y$ is reflexive it is enough to prove that 
$Y$ does not contain isomorphic copies of $\leins$ because by Rosenthal's theorem
\cite[Ch.\ XI]{Die-Seq}
in the absence of $\leins$ each bounded sequence
contains a weak Cauchy subsequence which then converges weakly by the
weak sequential completeness of $X$. But if $Y$ contained an isomorphic
copy of $\leins$ then by James' distortion theorem  it would contain also an
almost isometric copy of $\leins$.
By \cite[Cor.\ 3]{Pfi-Fixp} $Y$ would contain an asymptotic copy of $\leins$.
It's canonical normalized basis would $\tau_{\mu}$-converge to $0$;
finally, since $\tau_{\mu}$ and $\tau_{\norm{\cdot}}$
coincide on $\ball{Y}$
this basis would even converge in norm to $0$, a contradiction
which proves that $Y$ is reflexive.
\ebew\bigskip\\
Lemma \ref{lem3.1} is the technical key for the other results
of this section. It corresponds to \cite[IV.3.1]{HWW}.
\begin{lemma}\label{lem3.1}
Let $X$ be an L-embedded Banach space
with L-projection $P$ from $X''$ onto $X$.
Then for every net $(x_{\gamma})_{\gamma\in\Gamma}$ in $X$ \wst-converging
to $x''\in X''\backslash X$ there is a bounded sequence $(y_n)$ in 
$\mbox{co}\, \{x_{\gamma}\mdE \gamma\in\Gamma\}$ such that
the sequence $(y_n -Px'')$ spans $\leins$ asymtotically isometrically.
\end{lemma}
{\em Proof}:
Set $x=Px''$, $x_s=x''-Px''$. Choose a net
$(z_{\gamma})_{\gamma\in\Gamma}$ in $X$
such that $z_{\gamma}\genwst x_s$ and
$\norm{z_{\gamma}}=\norm{x_s}$.
We assume without loss of generality that both nets $(x_{\gamma})$, $(z_{\gamma})$
are indexed by the same directed set $\Gamma$.
Then $x_{\gamma}-z_{\gamma}\genw x$.\\
The idea of the proof is that on one hand by the theorem of
Hahn-Banach the net
$(x_{\gamma}-z_{\gamma})$ admits convex combinations which
converge to $x$
in norm and that on the other hand
by a slight modification of Godefroy's construction the corresponding 
convex combinations of the $x_{\gamma}$ can be chosen so to span
$\leins$ asymtotically. Here are the details.\\
Since $x''\not\in X$ we have $x_s\neq 0$ and thus
without loss of generality we suppose $\norm{x_s}=\norm{z_{\gamma}}=1$.
Let $(\delta_n)$ be a sequence of numbers in $]0,1[$
convergent to $0$.
Set $\eta_1=\frac{1}{4}\delta_1$ and
$\eta_{n+1}=\frac{1}{4}\min(\eta_n,\delta_{n+1})$ for
$n\in\N$.
By induction over $n\in \N$ we will construct
finite sets $A_n\subset \N$, finite sequences $(\lambda_k)_{k\in A_n}$
in $[0,1]$ and $(\gamma_k)_{k\in A_n}$ in $\Gamma$
such that
\bgl
\sum_{k\in A_{n}} \lambda_{k}     &=&1, \,\,\,\,
             G_i\cap G_n=\emptyset,\,\, G_n\subset \Gamma_n\,\,\,\falle i<n
                                                    \label{gl31_1}\\
&&\norm{(y_n -x)-\sum_{k\in A_n}\lambda_k z_{\gamma_k}}<\eta_n.
\label{gl31_5a}\\
\sum_{i=1}^{n}\betr{\alpha_i}(1+\eta_i)
&\geq&
\Norm{\sum_{i=1}^{n} \alpha_{i} (y_i -x)} 
   \geq \Bigl(\sum_{i=1}^{n}(1-\delta_i)\betr{\alpha_i}\Bigr)
  +\eta_n \sum_{i=1}^{n}\betr{\alpha_i}
  \falle \alpha_{i}\in\C
                                                    \label{gl31_3}
\egl
where
\bglst
y_n=\sum_{k\in A_{n}} \lambda_{k} x_{\gamma_k}, \,\,\,\,
G_n=\{\gamma_k\mdE k\in A_n \}, \,\,\,
\Gamma_n=\Gamma\backslash\bigcup_{i=1}^{n-1} G_i
\;\;(=\Gamma \mbox{ if }n=1).
\eglst
For the first induction step $n=1$ we choose $x'\in\ball{X'}$
such that $1=\norm{x_s}\geq\Ree x_s(x')>1-\delta_1 +2\eta_1$.
The \wst-convergence of $(z_{\gamma})$ to $x_s$ yields $\beta_1\in\Gamma$
such that\bgl
\Ree x'(z_{\gamma}) > 1- \delta_1 +2\eta_1 \falle \gamma\succeq\beta_1.
\label{gl_Go6a}
\egl
The net $((x_{\gamma}-x)-z_{\gamma})_{\gamma\succeq\beta_1}$
\w-converges to $0$ thus by the theorem of Hahn-Banach we find
a convex combination $y_1=\sum_{k\in A_1}\lambda_k x_{\gamma_k}$
such that
$$ \norm{(y_1-x)-\sum_{k\in A_1}\lambda_k z_{\gamma_k}}<\eta_1.$$
Thus (\ref{gl31_3}, $n=1$) follows from
$\norm{y_1-x}\leq\norm{\sum_{k\in A_1}\lambda_k z_{\gamma_k}}+\eta_1\leq1+\eta_1$ and from
\bglst
\norm{y_1-x}
&>&\norm{\sum_{k\in A_1}\lambda_k z_{\gamma_k}}-\eta_1\\
&\geq& -\eta_1 +\Ree \sum_{k\in A_1}\lambda_k x'(z_{\gamma_k}) \\
&>& -\eta_1 +\sum_{k\in A_1}\lambda_k(1-\delta_1 +2\eta_1)
=1-\delta_1 +\eta_1.
\eglst
For the induction step $n\mapsto n+1$ we suppose $A_i\subset \Gamma$,
$(\lambda_{k})_{k\in A_i} \subset [0,1]$,
$G_i\subset \Gamma$ to be
constructed according to \Ref{gl31_1} - \Ref{gl31_3} for
$i=1,\ldots,n$.

First we consider the \w-convergent net
$(x_{\gamma}-z_{\gamma})_{\gamma\in\Gamma_{n+1}}$.
By the theorem of Hahn-Banach we can choose a finite
set $A_{n+1} \subset \N$,
numbers $(\lambda_{k})_{k\in A_{n+1}}\subset [0,1]$ and indices
$(\gamma_k)_{k\in A_{n+1}}\subset\Gamma_{n+1}$ such that 
\Ref{gl31_1} and \Ref{gl31_5a} hold for $n+1$.
Together with (\ref{gl31_3}, $n$) this gives the first inequality of
(\ref{gl31_3}, $n+1$).

We fix an element $\alpha=(\alpha_{i})$ in the unit sphere of
$\leins_{n+1}$ such that $\alpha_{n+1}\neq 0$ and use the L-decomposition
of $X''=X\oplus_1 X_s$ in order to get
\bglst
\Norm{\Bigl(\sum_{i=1}^n \alpha_i (y_i -x)\Bigr) 
+ \alpha_{n+1}x_s}
&=&\Norm{\sum_{i=1}^n \alpha_i (y_i -x)} + \norm{\alpha_{n+1}x_s}\\
&\stackrel{\Ref{gl31_3}}{\geq}&
\Bigl(\sum_{i=1}^{n}(1-\delta_i)\betr{\alpha_i}\Bigr)
         +\eta_n \Bigl(\sum_{i=1}^{n}\betr{\alpha_i}\Bigr)
         +\betr{\alpha_{n+1}}\\
&=& \Bigl(\sum_{i=1}^{n+1}(1-\delta_i)\betr{\alpha_i}\Bigr)
        +\eta_n -(\eta_n -\delta_{n+1})\betr{\alpha_{n+1}}\\
&\geq& \Bigl(\sum_{i=1}^{n+1}(1-\delta_i)\betr{\alpha_i}\Bigr)
        +\min(\eta_n,\delta_{n+1})\\
&=& \Bigl(\sum_{i=1}^{n+1}(1-\delta_i)\betr{\alpha_i}\Bigr)
     +4\eta_{n+1}
\eglst
because $\norm{\alpha}=1$ and $\betr{\alpha_{n+1}}\leq1$.
Thus there is $x'\in\ball{X'}$ (depending on $\alpha$) such that
\bglst
\Norm{\Bigl(\sum_{i=1}^n \alpha_i (y_i -x)\Bigr) 
+ \alpha_{n+1}x_s}
&\geq&
\mbox{Re}{\Bigl(\Bigl(\sum_{i=1}^n \alpha_i (y_i -x)\Bigr) 
+ \alpha_{n+1}x_s\Bigr)(x')}\\
&>&\Bigl(\sum_{i=1}^{n+1}(1-\delta_i)\betr{\alpha_i}\Bigr)+3\eta_{n+1}.
\eglst
Then the \wst-convergence (along $\gamma\in\Gamma_{n+1}$) of
$\Bigl((\sum_{i=1}^n \alpha_i (y_i -x)) + \alpha_{n+1}z_{\gamma}\Bigr)$ to
$(\sum_{i=1}^n \alpha_i (y_i -x)) + \alpha_{n+1}x_s$ yields
$\beta\in\Gamma_{n+1}$ (depending on $\alpha$ and $x'$) such that
\bglst
\mbox{Re}\,\,\Bigl({x'\Bigl((\sum_{i=1}^n \alpha_i (y_i -x)) 
             + \alpha_{n+1}z_{\gamma}\Bigr)}\Bigr)
>
\Bigl(\sum_{i=1}^{n+1}(1-\delta_i)\betr{\alpha_i}\Bigr)+3\eta_{n+1}
 \,\,\falle\gamma\succeq\beta.
\eglst
Choose  a finite $\eta_{n+1}$-net
$(\alpha^{l})_{l=1}^{L_{n+1}}$ in the unit sphere
of $\leins_{n+1}$
in the sense that for each $\alpha$ in the unit sphere
of $\leins_{n+1}$ there is $l\leq L_{n+1}$ such that
$\norm{\alpha-\alpha^{l}} = \sum_{i=1}^{n+1}\betr{\alpha_i-\alpha_i^{l}}
<\eta_{n+1}$.
Then we may repeat the reasoning above finitely many times for
$l=1,\ldots, L_{n+1}$
in order to get $\beta_{n+1}\in\Gamma_{n+1}$
and $x_{l}'\in\ball{X'}$ for $l\leq L_{n+1}$ such that
\bgl
\mbox{Re}\left(x_l'\Bigl((\sum_{i=1}^{n}\alpha_i^l (y_i -x))
+\alpha_{n+1}^l z_{\gamma}\Bigr)\right)
> \Bigl(\sum_{i=1}^{n+1}(1-\delta_i)
\betr{\alpha_i^{l}}\Bigr)+3\eta_{n+1}
       \,\,\falle l\leq L_{n+1},\,\,\,\gamma\succeq\beta_{n+1}.
                             \label{gl31_4}  
\egl
For each $l\leq L_{n+1}$ we get that
\bgl
\Norm{\sum_{i=1}^{n+1}\alpha_i^l(y_i-x)}
&\stackrel{\Ref{gl31_5a}}{\geq}&
 -\eta_{n+1}    +\Norm{\Bigl(\sum_{i=1}^{n}\alpha_i^l(y_i-x)\Bigr)
              +\alpha_{n+1}^l \sum_{k\in A_{n+1}}\lambda_k z_{\gamma_k}}\\
&=& -\eta_{n+1} +   \Norm{\sum_{k\in A_{n+1}}\lambda_k
         \Bigl((\sum_{i=1}^{n}\alpha_i^l(y_i-x))+
                                    \alpha_{n+1}^l z_{\gamma_k}\Bigr)}\\
&\geq& -\eta_{n+1} 
     +\sum_{k\in A_{n+1}} \lambda_k \Ree x_l'\Bigl((\sum_{i=1}^{n}\alpha_i^l(y_i-x))+
                                    \alpha_{n+1}^l z_{\gamma_k}\Bigr)\\
&\stackrel{\Ref{gl31_4}}{\geq}&
\Bigl(\sum_{i=1}^{n+1}(1-\delta_i)
\betr{\alpha_i^{l}}\Bigr)+2\eta_{n+1} \label{gl31_5}
\egl
For an arbitrary $\alpha$ in the unit sphere
of $\leins_{n+1}$ choose $l\leq L_{n+1}$ such that
$\norm{\alpha-\alpha^{l}}<\eta_{n+1}$.
Then
\bgl
\Norm{\sum_{i=1}^{n+1}\alpha_i (y_i -x)}
&\geq&
\Norm{\sum_{i=1}^{n+1}\alpha_i^{l} (y_i -x)} -
   \Norm{\sum_{i=1}^{n+1}(\alpha_i - \alpha_i^{l}) (y_i -x)}
\nonumber\\
&\geq&
\Bigl(\sum_{i=1}^{n+1}(1-\delta_i)\betr{\alpha_i}\Bigr)+2\eta_{n+1}
-\norm{\alpha-\alpha^{l}}
\nonumber\\
&\geq&
\Bigl(\sum_{i=1}^{n+1}(1-\delta_i)\betr{\alpha_i}\Bigr)+\eta_{n+1}
\nonumber\\
&=&
\Bigl(\sum_{i=1}^{n+1}(1-\delta_i)\betr{\alpha_i}\Bigr)
      +\eta_{n+1}\sum_{i=1}^{n+1}\betr{\alpha_i}.
\nonumber
\egl
This extends to all scalars $\alpha_i\in\C$ and thus
ends the induction.
The sequence $(y_n)$ is bounded because of \Ref{gl31_5a}.
By \Ref{gl31_3} the sequence $(y_n -x)$ is easily seen to
span $\leins$ asymptotically.
This ends the proof.\ebew\bigskip\\
Remark: The proof yields not only
$y\in\mbox{co}\, \{x_{\gamma}\mdE \gamma\in\Gamma\}$
but separated blocks
$y_n=\sum_{k\in A_n}\lambda_k x_{\gamma_k}$ where the sets
$\{x_{\gamma_k}\mdE k\in A_n\}$ are pairwise disjoint.
Moreover  one can obtain, given a sequence $(\gamma'_n)$ in $\Gamma$,
that $x_{\gamma_k}\succeq \gamma'_n$ for $k\in A_n$.
\bigskip\\
In general the unit ball of $X$ is not $\tau_{\mu}$-compact;
the Rademacher functions $r_n$ in $\Leins([0,1])$ which are bounded
without having a measure convergent subsequence provide
a counterexample. [If $(r_n)$ contained
a measure convergent subsequence this subsequence would
admit a norm convergent subsequence
by Proposition \ref{prop_reflex} (a) because
$(r_n)$ spans $l^2$ isomorphically.] But we have:
\begin{coro}\label{coro_compakt}
Every bounded sequence in an L-embedded Banach space admits
a sequence of convex combinations which converges with respect to
the measure topology.
\end{coro}
{\em Proof}:
Let $(x_n)$ be a bounded sequence in an L-embedded space $X$ and let
$(x_{n_{\gamma}})$ be a universal net that \wst-converges to
$x''\in X''$ by the \wst-compactness of $\ball{X''}$.
If $x''\in X$ then this net
$w$-converges, admits norm-convergent convex combinations and
we are done in this case.
Otherwise $x''$ lies in $X''\backslash X$ and one applies
Lemma \ref{lem3.1} to get
a sequence  $(y_n)$ in $\mbox{co}\, \{x_n\mdE n\in\N\}$ which
$\tau_{\mu}$-converges to $Px''$.
\ebew\medskip\\
Corollary \ref{coro_compakt} corresponds to \cite[p.\ 202]{HWW}.
In this context there is a natural\\
{\bf Question}: Does Komlos' theorem hold accordingly?
More precisely, given
a bounded sequence in an L-embedded space, does it admit a aubsequence
whose Cesaro (=arithmetic) means converge with respect to the measure
topology?\\
Note that by (a) of Proposition \ref{prop_reflex}
Komlos' theorem implies the weak Banach-Saks property (which
by definition claims that a \w-convergent sequence admits a subsequence
whose Cesaro means converge in norm, see for example
\cite{Beauz}, \cite[p.\ 112, 121]{Die-Seq}, \cite{BelDie}).
By Rosenthal's $\leins$-theorem the weak Banach-Saks property is
also half a converse to Komlos' theorem, that is by
Rosenthal's $\leins$-theorem a bounded sequence in an L-embedded space
admits a subsequence which is either equivalent to the standard basis
of $\leins$ or converges weakly; but if one supposes the 
weak Banach-Saks property to hold then in the second case of a
\w-convergent sequence there are Cesaro means that converge in norm
whence with respect to the measure topology.\\
There is another related\\
{\bf Question}: Does the Kadec-\Pel\ subsequence
decomposition (sometimes also called the Kadec-\Pel\
splitting lemma) hold accordingly?
This lemma says that a bounded sequence $(f_n)$ in $\Leins([0,1])$ admits a
subsequence $(f_{n_k})$ which can be decomposed in the following sense:
there are two bounded sequences $(g_k)$, $(h_k)$ in $\Leins([0,1])$ such that
$f_{n_k}=g_k + h_k$, the $g_k$ are pairwise disjoint, and $(h_k)$ converges
weakly.
In a recent preprint Randrianantoanina \cite{Randri-Kad-Pel}
showed the Kadec-\Pel\ subsequence
decomposition for preduals of von Neumann algebras.
Since the latter are known to have the weak Banach-Saks property
\cite{BelDie}, Komlos' theorem follows almost immediately from
Randrianantoanina's result for von Neumann preduals, see Proposition
\ref{theoKomlosvN} below.
\bigskip\\
The following theorem generalizes a theorem of Buhvalov-Lozanovskii
(\cite{Buh-Lo}, \cite[IV.3.4]{HWW}).
As in \cite[IV.3.4]{HWW} the implication (ii)$\folgt$(i)
holds also for unbounded $C$.
\begin{satz}\label{theoBL}
Let $X$ be an L-embedded Banach space with L-projection $P$
(on $X''$ with range $X$) and
endowed with its abstract measure topology $\tau_{\mu}$.
For a norm closed bounded convex set $C\subset X$ the following two
assertions are equivalent.
\begin{description}
\item{\makebox[1.8em]{(i)}}
$P\,\overline{C}^{w^{*}}=C$
where \wst\ refers to the \wst-topology of $X''$.
\item{\makebox[1.8em]{(ii)}}
$C$ is $\tau_{\mu}$-closed.
\end{description}
\end{satz}
{\em Proof}:
(i)$\folgt$(ii) Take $c_n\in C$, $x\in X$ with $c_n\gentaumu x$.
It is enough to show that $x\in C$ because closedness and
sequential closedness coincide for $\tau_{\mu}$.
If the $\tau_{\mu}$- null sequence $(c_n -x)$
contains a norm convergent subsequence then we are done
because in the norm topology $C$ is closed.
Otherwise an appropriate subsequence $(c_{n_k} -x)$ spans $\leins$
asymptotically and $\inf \norm{c_{n_k}-x}>0$.
By \wst-compactness of $\overline{C}^{w^{*}}$ there is a net
$(c_{n_{k_{\gamma}}}-x)_{\gamma\in\Gamma}$ on
${\cal S}=\{c_{n_k}-x\mdE k\in\N\}$
that \wst-converges to $\overline{c}-x\in\overline{C}^{w^{*}}-x$.
By Lemma \ref{lem_HP} we have that
$\overline{c}-x\in X_s$ because $C$ is bounded.
Thus $P(\overline{c}-x)=0$ and $x=Px=P\overline{c}\in C$
by hypothesis.\\
(ii)$\folgt$(i)
This implication is essentially Lemma \ref{lem3.1}:
If $(c_{\gamma})_{\gamma\in\Gamma}$
is a \wst-convergent net in $C$ with limit $\overline{c}$ it is enough
to prove that $P\overline{c}\in C$ because the inclusion
$C\subset P\overline{C}^{w^{*}}$ is trivial.
If $\overline{c}\in X$ then there is a sequence of convex combinations
of the $c_{\gamma}$ convergng to $\overline{c}$ in norm
whence $\overline{c}\in C$.
Otherwise, if $\overline{c}\in X''\backslash X$,
by Lemma \ref{lem3.1} there is a sequence
$(d_n)\in \,\mbox{co}\,\{c_{\gamma}\mdE \gamma\in\Gamma\}\subset C$ which
$\tau_{\mu}$-converges to $P\overline{c}$ hence $P\overline{c} \in C$ because $C$ is $\tau_{\mu}$-closed.
\ebew\bigskip\\
Corollary \ref{coroBL}
corresponds to \cite[IV.3.5]{HWW}) for
the case $X=\Leins(\Omega,\Sigma,\mu)$.
\begin{coro}\label{coroBL}
Let $X$ be an L-embedded Banach space endowed with its abstract
measure topology $\tau_{\mu}$.
Then a norm closed subspace $Y\subset X$ is L-embedded if and only if its
unit ball $\ball{Y}$ is $\tau_{\mu}$-closed.
\end{coro}
The proof is immediate from Theorem \ref{theoBL} with $C=\ball{Y}$
and from Li's criterion Lemma  \ref{lem_Li}.
\ebew\bigskip\\
Let $X$ be a Banach space admitting an abstract measure topology
$\tau_{\mu}$. Then we define
$$X^{\#}=\{x'\in X'\mdE x'\eing{\ball{X}}
\,\,\mbox{is}\,\,\tau_{\mu}\mbox{continuous}\}.$$
Remarks:\\
1. $X^{\#}$ is a closed subspace of $X'$. (The proof is left to the reader.)\\
2. If $X$ is a subspace of an L-embedded Banach space
then one has $X^{\#}=X'$
if and only if $X$ does not contain copies of $\leins$.
For, in the absence of $\leins$, $\tau_{\mu}$ coincides with the norm topology
hence $X^{\#}=X'$.
Conversely, if $X$ contains a copy of $\leins$ then
by James' destortion theorem it contains
also an almost isometric copy $U$ of $\leins$ spanned by a normalized 
basis $(u_n)$.
Since $X$ is contained in an L-embedded space, by \cite[Cor.\ 3]{Pfi-Fixp}
we (may pass to an appropriate subsequence and) suppose that
$u_n\gentaumu0$.
Let $x'\in X'$ be a Hahn-Banach extension of the functional on
$U$ which corresponds to $\eins\in\lunendl$. Then $x'$ is not
$\tau_{\mu}$-continuous on $\ball{X}$ since $u_n\gentaumu0$
but $x'(u_n)\gen1$. (Compare also with \cite[Rem.\ (b)\ p.\ 186]{HWW}.)
\ebew
\begin{proposition}\label{prop3.9}
Let $X$ be an L-embedded Banach space (with L-decomposition
$X''=X\oplus_1 X_s$) endowed with its abstract
measure topology $\tau_{\mu}$. Then
\bglst
X^{\#} = (X_s)_{\bot} \mbox{ (}\subset X'\mbox{)}.
\eglst
\end{proposition}
{\em Proof}:
"$\subset$"
Take $x'\in X^{\#} $ and $x_s \in X_s$.
To prove the inclusion we show that $x_s(x')=0$.

Let $(x_{\gamma})$ be a net that \wst-converges to $x_s$
with $\norm{x_{\gamma}}=\norm{x_s}$.
But then, by Lemma \ref{lem_HP}
there is a sequence $(x_{\gamma_n})$
that spans $\leins$ asymptotically.
Hence $x_{\gamma_n}\gentaumu 0$ and
$x'(x_{\gamma_n})\gen0$ since $x'$ is $\tau_{\mu}$-continuous on
bounded sets. This proves that $x''(x')=0$.
(In passing we note that $x_{\gamma}\gentaumu 0$
by (h) of Lemma \ref{lem_elementar} but that one can
not infer from this that $x''(x')=0$ because
it is not clear whether a
$\tau_{\mu}$-convergent net has a unique limit.)\\
"$\supset$"
Assume that there is $x'\in \bigl(X_s)\bigr)_{\bot}$
that is not $\tau_{\mu}$-continuous on $\ball{X}$.
Then by the definiton of $\tau_{\mu}$ there are $\eps>0$ and a sequence
$(x_n)$ in $\ball{X}$ such that $x_n\gentaumu 0$
but $\betr{x'(x_n)}>\eps$ for all $n\in\N$.
Still by definiton of $\tau_{\mu}$ and because $x'$ is norm-continuous
we suppose that $(x_n)$ spans $\leins$ almost isometrically.
By \wst-compactness of $\ball{X''}$
there exists a \wst-accumulation point $x''\in \ball{X''}$ of
$\{x_n\mdE n\in\N\}$. Let $(x_{n_\gamma})$ be a net \wst-converging
to $x''$. By Lemma \ref{lem_HP},
$x''\in X_s$.
Thus $x'(x_{n_{\gamma}})\gen x''(x')=0$ by hypothesis.
This contradicts $\betr{x'(x_n)}>\eps$ and proves that $x'$ is
$\tau_{\mu}$-continuous on $\ball{X}$.
\ebew\bigskip\\
Theorem \ref{theo3.10} (see also \cite[IV.3.10]{HWW})
was proved in \cite{G-Li} for nicely placed (=L-embedded) subspaces
$X$ of $\Leins(\Omega,\Sigma,\mu)$, $\mu$ finite.
For its proof Proposition \ref{prop3.9} plays the same r\^ole as \cite[IV.3.9]{HWW}
in the proof of \cite[IV.3.10]{HWW}.
Recall that if
an L-embedded space admits a predual then this predual need not
be M-embedded (\cite[p.\ 102]{HWW}).
\begin{satz}\label{theo3.10}
Let $X$ be an L-embedded Banach space endowed with its abstract
measure topology $\tau_{\mu}$. The following assertions are
equivalent.
\begin{description}
\item{\makebox[1.8em]{(i)}}
$X$ is (isometrically isomorphic to) the dual of an M-embedded
Banach space $Z$.
\item{\makebox[1.8em]{(ii)}}
$X^{\#}$ separates $X$.
\end{description}
If (i) and (ii) are satisfied $Z$ is (isometrically isomorphic to)
$X^{\#}$.
\end{satz}
{\em Proof}:
(i)$\folgt$(ii):
By \cite[III.1.3]{HWW} the L-projections on $X''$ and on $Z'''$ with kernel
$X_s=Z^{\bot}\subset X''$ can be identified.
Thus $Z=(Z^{\bot})_{\bot}=(X_s)_{\bot}=X^{\#}$ by
Proposition \ref{prop3.9}; in particular, $X^{\#}$ separates $X$.\\
(ii)$\folgt$(i):
By \cite[IV.1.9]{HWW} it is enough to show that $X_s$ is \wst-closed in $X''$;
then an M-embedded predual of $X$ exists and is isometrically isomorphic to
$(X_s)_{\bot}$ whence to $X^{\#}$ by Proposition \ref{prop3.9}.
To see that $X_s$ is \wst-closed we take an element
\bglst
\overline{x}=x+ x_s \in \overline{X_s}^{w^*} =((X_s)_{\bot})^{\bot}
=(X^{\#})^{\bot}
\eglst
with $x\in X$, $x_s \in X_s$.
Then $x=\overline{x} - x_s \in (X^{\#})^{\bot} \cap X =(X^{\#})_{\bot}=\{0\}$
where the latter  equality comes from the fact that
$X^{\#}$ separates $X$.
Thus $\overline{x}=x_s \in X_s$ which proves that
$X_s$ is \wst-closed in $X''$.
\ebew\medskip\\
Analogously to the $\leins$-case we say that a sequence $(x_n)$ of nonzero
elements in a Banach space $X$ spans $c_0$ almost (respectivley asymptotically)
isometrically if  there exists a sequence
$(\delta_m)$ in $[0,1[$ tending to $0$ such that 
$(1-\delta_m)\max_{m\leq n\leq m'}\betr{\alpha_n}
\leq
\Norm{\sum_{n=m}^{m'} \alpha_n \frac{x_n}{\norm{x_n}}}
\leq
(1+\delta_m)\max_{m\leq n\leq m'}\betr{\alpha_n}$
for all $m\leq m'$
(respectively such that
$\max_{n\leq m}(1-\delta_n)\betr{\alpha_n}
\leq \Norm{\sum_{n=1}^{m}\alpha_n \frac{x_n}{\norm{x_n}}}
\leq \max_{n\leq m}(1+\delta_n)\betr{\alpha_n}$ for all $m\in\N$).

It follows from the proof of \cite[Th.\ 2]{DJLT} that the dual of an
asymptotic $c_0$-copy is an asymptotic $\leins$-copy.
A similar argument shows that this remains true if "asymptotic" is replaced by
"almost isometric".
Analogously we get
\begin{lemma}\label{theo_dual}
Let $X=[x_n]$ be an almost isometric copy of $\leins$ spanned by the
normalized basis $(x_n)$.
Let $x_n' \in X'$ be the biorthogonal functionals
(that is $x_n'(x_k)=\delta_{n,k}$).
Then $(x_n')$ spans $c_0$ almost isometrically.\\
The same holds with "almost" replaced by "asymptotically".
\end{lemma}
{\em Proof}: First we deal with the case where $x_n\almleins$.\\
By hypothesis there is a null sequence $(\delta_m)\subset[0,1[$ sucht that
\bglst
(1-\delta_m)\sum_{m}^{\infty}\betr{\beta_n}
\leq 
\Norm{\sum_{m}^{\infty} \beta_n x_n}
\leq
\sum_{m}^{\infty}\betr{\beta_n}
\eglst
for all scalars $\beta_n$. Let $m\leq m'$ be arbitrary in $\N$.
Since $\norm{x_n}=1$
we have the first inequality of
\bgl
\max_{m\leq n\leq m'} \betr{\alpha_n}\leq 
\Norm{\sum_{n=m}^{m'} \alpha_n x_n'}
\leq
(1-\delta_m)^{-1} \max_{m\leq n\leq m'}\betr{\alpha_n}
\label{gl_dual2}
\egl
for all scalars $\alpha_n$.
For the second inequality we take any $x\in \ball{X}$ of the form
$x=\sum \beta_n x_n$; then 
\bglst
\Betr{\sum_{n=m}^{m'}\alpha_n x_n'(x)}
&=& \Betr{\sum_{n=m}^{m'}\alpha_n \beta_n}
\leq \max_{m\leq n\leq m'}\betr{\alpha_n}
         \sum_{n=m}^{m'}\betr{\beta_n}\\
&\leq& (1-\delta_m)^{-1} \max_{m\leq n\leq m'}\betr{\alpha_n}
\eglst
whence the second inequality of \Ref{gl_dual2}.
It follows from \Ref{gl_dual2} that $(x_n')$ spans $c_0$ almost isometrically.\\
The case in which $x_n\asyleins$ is proved similarly.
\ebew\medskip\\
\begin{proposition}\label{prop_leins_dual}
An L-embedded almost isometric copy of $\leins$ is the dual of
an M-embedded space which is almost isometric to $c_0$.\\
The statement remains true when "almost" is replaced by "asymptotically".
\end{proposition}
{\em Proof}:
Let $X$ be an L-embedded almost isometric $\leins$-copy with a
normalized canonical basis $(x_n)$.
Let $(x_n')$ be the biorthogonal functionals of $(x_n)$ that is
$x_n'(x_m)=\delta_{n,m}$ for $n,m\in\N$.\medskip\\
{\bf Sublemma}
{\em $x_n' \in X^{\#}$ for all $n\in\N$.}\smallskip\\
{\em Proof} of the Sublemma:
Suppose there is $n_0 \in\N$ such that $x'=x_{n_0}\not\in X^{\#}$.
Then there is a sequence $(y_n)\subset \ball{X}$ and there is $\eps>0$
such that 
\bglst
\norm{y_n}=1,\;\;\;
y_n\asyleins, \;\;\;
\betr{x'(y_n)}>\eps
\mbox{ for all } n\in\N.
\eglst
By \cite[Lem.\ 1]{Pfi-Fixp} $Y=[y_n]$ is L-embedded and $Y^{\bot\bot}
=Y\oplus_1 Y_s$ with $Y_s=Y^{\bot\bot}\cap X_s$.
Let $\delta>0$ be arbitrary for the moment.
Set $X_0=[x_n]_{n\geq m_0}$ where $m_0$ is such that 
$(x_n)_{n\geq m_0}\isoleins{1-\delta}$. $X_0$ is L-embedded by
\cite[Lem.\ 1]{Pfi-Fixp} and we write
$X_0^{\bot\bot}=X_0 \oplus_1 (X_0)_s$ with $(X_0)_s=X_0^{\bot\bot}\cap X_s$.
Let $y_s\in X''$ be a \wst-accumulation point of $\{y_n \mdE n\in\N\}$.
By Lemma \ref{lem_HP} we have $\norm{y_s}=1$ and $y_s\in Y_s$.
We have $y_s\in (X_0)_s$ because $X_0$ is co-finite-dimensional in $X$;
furthermore, $\betr{y_s(x')}\geq \eps$.
Let $(z_{\gamma})\subset \ball{X_0}$ be a normalized net \wst-converging to 
$y_s$. After passing to an appropriate subnet we  suppose that
$\betr{(y_s-z_{\gamma})(x')}<\eps/2$ for all $\gamma$.
By Lemma \ref{lem_HP} one can extract a sequence $(z_{\gamma_n})$ that spans
$\leins$ asymptotically.
Then $\betr{x'(z_{\gamma_n})}\geq\eps/2$ for all $n$.
We define an isomorphism $T:X_0\gen [e_n]_{n\geq m_0}$ by $x_n\mapsto e_n$
where $(e_n)$ is the standard basis of $\leins$.
Then $\norm{T}\leq (1-\delta)^{-1}$ and $\norm{T^{-1}}\leq1$.
There is $m_1\in\N$ such that
$(z_{\gamma_n})_{n\geq m_1}\isoleins{1-\delta}$
since $z_{\gamma_n}\almleins$.
Hence, with the notation $f_n=Tz_{\gamma_n}$, we get
$(f_n)_{n\geq m_1}\isoleins{1-2\delta}$ because
\bglst
(1-2\delta)\sum_{n=m_1}^{\infty}\betr{\alpha_n}
&<&(1-\delta)^2 \sum_{n=m_1}^{\infty}\betr{\alpha_n}
\leq (1-\delta)\sum_{n=m_1}^{\infty}\frac{\betr{\alpha_n}}{\norm{f_n}}\\
&\leq&
\Norm{\sum_{n=m_1}^{\infty}\frac{\alpha_n}{\norm{f_n}}z_{\gamma_n}}
\leq
\Norm{\sum_{n=m_1}^{\infty}\alpha_n\frac{f_n}{\norm{f_n}}}.
\eglst
Then by \cite[Th.\ B]{Dor}
or \cite[L.\ 10, 6]{Pfi-vN} there is
a sequence $(\tilde{f}_k)\subset\leins$ of pairwise disjointly supported
elements of $\leins$ and a subsequence $(f_{n_k})$ such that
\bglst
\norm{f_{n_k}-\tilde{f}_k} <\delta'
\eglst
where $\delta'\gen0$ as $\delta\gen0$.
Set $\tilde{z}_k=T^{-1}\tilde{f}_k$.
Then
\bglst
\norm{\tilde{z}_k-z_{\gamma_{n_k}}}=
\norm{T^{-1}\tilde{f}_k -T^{-1}f_{n_k}}
\leq \norm{\tilde{f}_k -f_{n_k}}<\delta'.
\eglst
Now we choose $\delta>0$ small enough in order to have $\delta'\norm{x'}<\eps/4$.
Hence
\bglst
\betr{x'(\tilde{z}_k)}\geq \frac{\eps}{2}-
\betr{x'(\tilde{z}_k -z_{\gamma_{n_k})}}
\geq \frac{\eps}{2}-\delta'\norm{x'}\geq \frac{\eps}{4}
\eglst
for all $k\in\N$.
On the other hand, for $e'=e_{n_0}'=(T^{-1})'(x')$
we have $e'(e_n)=\delta_{n_0, n}$ that is 
\bglst
x'(\tilde{z}_k )=e'(\tilde{f}_k)=0
\eglst
for all but possibly one $k\in\N$
because the $\tilde{f}_k$ are pairwise disjoint.
This contradiction proves the Sublemma.
\ebew\smallskip\\
Since the biorthogonal functionals separate $X$ the  Sublemma and Theorem
\ref{theo3.10} imply that $X^{\#}$ is M-embedded and $X=(X^{\#})'$.
The Sublemma states that $[x_n']\subset X^{\#}$. In fact one has
$[x_n']= X^{\#}$. To see this let $z'\in X^{\#}\backslash[x_n']$.
Since $X'$ is isomorphic to $\lunendl$ there is an infinite set $N'\subset\N$
and there is $\eps>0$ such that $\betr{z'(x_n)}>\eps$ for all $n\in N'$.
But $(x_n)_{n\in N'}\almleins$ whence $(x_n)_{n\in N}\asyleins$ for an appropriate
infinite set $N\subset N'$ by \cite[Cor.\ 3]{Pfi-Fixp}.
This means that $(x_n)_{n\in N}$ $\tau_{\mu}$-converges to $0$ hence $z'\not\in X^{\#}$.

We have proved that if $X$ is L-embedded and almost isometric to $\leins$ then
$X$ is the dual of the M-embedded space $X^{\#}$ and $X^{\#}=[x_n']$.
To end the proof it just remains to apply Proposition \ref{theo_dual}.\\
The arguments are similar for the case in which $x_n\asyleins$.
\ebew\bigskip\bigskip\\
{\bf\S 5 On some results of Godefroy, Kalton, Li}\\
First we deal with \cite[Prop.\ 2.1]{GoKaLi}. The first part of that proposition
says that $X^{\#}$ where $X$ is a nicely placed (=L-embedded) subspace of
$\Leins(\mu)$ is always M-embedded, not only in the situation
of Theorem \ref{theo3.10}.
\begin{proposition}\label{theo2.1a}
Let $X$ be an L-embedded Banach space.
Then $X^{\#}$ is M-embedded.
\end{proposition}
{\em Proof}:
With the usual notation $X''=X\oplus_1 X_s$ and with Proposition \ref{prop3.9}
we have $$X^{\#}=(X_s)_{\bot}\subset X'.$$
We set $$Z=\overline{X_s}^{w^*}=(X^{\#})^{\bot}\subset X''$$
and $$Y=X\cap \overline{X_s}^{w^*}=(X^{\#})_{\bot}\subset X.$$
By Corollary \ref{coroBL}, $Y$ is L-embedded because its unit ball
\bglst \ball{Y}=\ball{X}\cap\bigcap_{x'\in X^{\#}}\mbox{ker}\,x'\eglst
is $\tau_{\mu}$-closed. Hence $Y^{\bot\bot}=Y\oplus_1 (Y^{\bot\bot}\cap X_s)$
by Lemma \ref{lem_Li}.
Now the fact $Z=Y\oplus_1 X_s$ and the fact that $X^{\#}$ is an M-ideal
in its bidual $Y^{\bot}$ can be deduced exactly as in the proof of
\cite[Prop.\ 2.1]{GoKaLi}.\ebew\bigskip\\
For the proof of Proposition \ref{theo2.1b} we recall property $(m_1^*)$.
In \cite{KaWe} a separable Banach space $Z$ is defined
to have property $(m_1^*)$ if for all $z', z_n'\in Z'$
\bgl \limsup\norm{z'+z_n'}=\norm{z'}+\limsup\norm{z_n'}\label{glBem1}\egl
whenever $z_n'\genwst0$.
Analogously a separable Banach space $X$ is defined to have  property
$(m_1)$ if for all $x, x_n \in X$
\bgl \limsup\norm{x+x_n}=\norm{x}+\limsup\norm{x_n}\label{glBem2}\egl
whenever $x_n\genw 0$.\medskip\\
The second part of \cite[Prop.\ 2.1]{GoKaLi} reads as follows in our context.
\begin{proposition}\label{theo2.1b}
Let $X$ be a separable L-embedded Banach space.
If $\ball{X}$ is $\tau_{\mu}$-sequentially compact then for any
$\eps>0$ there is a
subspace $X_{\eps}$ of $c_0$ such that $\dist(X^{\#},X_{\eps})<1+\eps$.
\end{proposition}
{\em Proof}:
Exactly as in \cite{GoKaLi} one distinguishes three steps: Firstly one proves that $X^{\#}$ has
property $(m_1^*)$, secondly one deduces from this property $(m_{\infty}^*)$
and thirdly it remains to apply \cite[Th.\ 3.5]{KaWe}.
Only the first step must be modified a bit.\\
From the proof of Proposition \ref{theo2.1a} we know that
$(X^{\#})'=X/Y$ where $Y=X\cap \overline{X_s}^{w^*}=(X^{\#})_{\bot} \subset X$.
Let $(u_n)\subset X/Y=(X^{\#})'$ be a \wst-null sequence. We denote by $Q:X\gen X/Y$
the quotient map.
Let $(x_n)\subset X$ be a bounded sequence such that $Qx_n=u_n$. By hypothesis there is a
$\tau_{\mu}$-convergent subsequence - still denoted by $(x_n)$ - such that
$x_n-x_0\asyleins$ where $x_0=\tau_{\mu}-\lim x_n$ and such that
$\lim\norm{x_0 -x_n}$ exists.
We have $x_0\in Y$ because for any $x'\in X^{\#}$ one has
\bglst x'(x_0)=\lim x'(x_m)=\lim x'(u_m)=0.\eglst
We have furthermore that
\bgl
\limsup\norm{y+x+(x_n-x_0)}
\geq \norm{y+x} + \lim\norm{x_n-x_0}\label{gl2.1b1}\egl
for all $x\in X$, $y\in Y$.
To see this, recall that $X$ is L-embedded, and that by Lemma \ref{lem_HP}
 each universal
net $(x_{n_\gamma}-x_0)$ \wst-converges to a limit $x_s\in X_s$ such that
$\norm{x_s}=\lim_{\gamma}\norm{x_{n_\gamma}-x_0}$ and
\bglst
\lim_{\gamma}\norm{y+x+(x_{n_\gamma}-x_0)}
\geq \norm{y+x+x_s}=\norm{y+x}+\norm{x_s}
\eglst
by \wst-continuity of the norm whence \Ref{gl2.1b1}.

Since $\norm{Qx}=\inf_{y\in Y}\norm{y+x}$ we deduce from \Ref{gl2.1b1} that 
\bglst
\limsup\norm{y+x+(x_n-x_0)}\geq \norm{Qx}+\lim\norm{u_n}
\eglst
and
\bglst
\limsup\norm{Qx+u_n}\geq \norm{Qx}+\lim\norm{u_n}
\eglst
which proves that $X^{\#}$ has property $(m_1^*)$.
The deduction of property $(m_{\infty}^*)$ and the conclusion via
\cite[Th.\ 3.5]{KaWe} do not depend on the measure topology and coincide therefore
with the arguments in \cite{GoKaLi}.\ebew\bigskip\\
The following remark gives a characterization of property $(m_1^*)$ in
L-embedded Banach spaces. We will not need it in the sequel and state it
only because the way we prove it by constructing asymptotic $\leins$-sequences
fits naturally in the main theme of this paper.
Note that the implication (i)$\folgt$(ii) holds for arbitrary
Banach spaces $Z$, that the implications
(i)$\folgt$(ii)$\gdw$(iii) hold whenever $Z$ is such that
its dual admits an abstract measure topology, and that the implication
(iii)$\folgt$(iv) does not need the M-embeddedness of $Z$.
Note furthermore that in Remark \ref{Bem} the separation assumption on $Z$
could be omitted because the definition of properties $(m_1^*)$ and $(m_1)$
makes sense also for non-separable spaces.\\
\begin{remark}\label{Bem}\mbox{}\\
(a) Let $Z$ be a Banach space such that its dual is L-embedded.
Then the following assertions are equivalent:
\begin{description}
\item{\makebox[1.8em]{(i)}}
$Z$ has property $(m_1^*)$.
\item{\makebox[1.8em]{(ii)}}
Each \wst-null sequence in $Z'$ admits of  a subsequence that converges to $0$
in norm or spans $\leins$ asymptotically.
\item{\makebox[1.8em]{(iii)}}
Each \wst-null sequence in $Z'$ is $\tau_{\mu}$-null.
\end{description}
If $Z$ is even separable and M-embedded then the assertions above are equivalent to
\begin{description}
\item{\makebox[1.8em]{(iv)}}
$\ball{Z'}$ is $\tau_{\mu}$-sequentially compact.
\end{description}\smallskip
(b)
An arbitrary Banach space $X$ has property $(m_1)$ if and only if
it has the Schur property.
\end{remark}
{\em Proof}: (a) We set $X=Z'$.
Sketch of\\
(i)$\folgt$(ii):
The proof ressembles the one of \cite[Th.\ 2]{Pfi-Fixp} the only
difference being that \Ref{glBem1} replaces the \wst-lower semicontinuity of
the norm and the L-embeddedness of $Z'$.
Let $(x_m)\subset X$ be a \wst-null sequence, suppose
without loss of generality that $\norm{x_m}=1$.
Let $(\delta_n)$ be a sequence of strictly positive numbers
converging to $0$.
Set $\eta_1=\frac{1}{6}\delta_1$ and
$\eta_{n+1}=\frac{1}{6}\min(\eta_n,\delta_{n+1})$ for
$n\in\N$.
By induction over $n\in\N$ one constructs $m_n\in\Gamma$
such that 
\bgl
\Bigl(\sum_{i=1}^{n}(1-\delta_i)\betr{\alpha_i}\Bigr)
+\eta_n \sum_{i=1}^{n}\betr{\alpha_i}
\leq
\norm{\sum_{i=1}^{n}\alpha_i x_{m_i}}
\leq         \sum_{i=1}^{n}\betr{\alpha_i}
\falle n\in\N, \,\,\alpha_i\in\C.
\label{glGo1}
\egl
For the induction step $n\mapsto n+1$
fix an element $\alpha=(\alpha_i)_{i=1}^{n+1}$ in the
unit sphere of $\leins_{n+1}$ such that $\alpha_{n+1}\neq 0$.
Then \Ref{glBem1} yields
\bglst
\lefteqn{\Norm{\Bigl(\sum_{i=1}^n \alpha_i x_{m_i}\Bigr) 
+ \alpha_{n+1}x_m} +\frac{1}{2}\min(\eta_n,\delta_{n+1})}\\
&\geq&
\Norm{\sum_{i=1}^n \alpha_i x_{m_i}} + \betr{\alpha_{n+1}}\\
&\stackrel{\Ref{glGo1}}{\geq}&
\Bigl(\sum_{i=1}^{n}(1-\delta_i)\betr{\alpha_i}\Bigr)
         +\eta_n \Bigl(\sum_{i=1}^{n}\betr{\alpha_i}\Bigr)
         +\betr{\alpha_{n+1}}\\
&\geq& \Bigl(\sum_{i=1}^{n+1}(1-\delta_i)\betr{\alpha_i}\Bigr)
        +\min(\eta_n,\delta_{n+1})
\eglst
whence
\bglst
\Norm{\Bigl(\sum_{i=1}^n \alpha_i y_{m_i}^{(i)}\Bigr) 
+ \alpha_{n+1}y_m^{(n+1)}}
\geq
\Bigl(\sum_{i=1}^{n+1}(1-\delta_i)\betr{\alpha_i}\Bigr) + 3\eta_{n+1}
\eglst
for infinitely many $m$.
This gives \Ref{glGo1}; the details are the same as in the proof
of \cite[Th.\ 2]{Pfi-Fixp} or of Lemma \ref{lem_Frechet}.\\
Note for the proof of part (b) below that we used the \wst-convergence of
$(x_n)$ only in order to apply property $(m_1^*)$ but not for the construction
of the $\leins$-basis itself.\\
(ii)$\folgt$(i):
Suppose that $Z$ satisfies (ii) without having $(m_1^*)$ and suppose that
$Z'$ is L-embedded.
Then there are $x\in X$, $\eps>0$ and a \wst-null
sequence $(x_n)\subset X$
such that $\lim\norm{x+x_n}$ and $\lim\norm{x_n}$ exist and
\bgl \eps+\lim\norm{x+x_n}< \norm{x}+\lim\norm{x_n}.\label{glBem3}\egl
Since \Ref{glBem3} excludes the case $\norm{x_n}\gen0$ we have - after passing
to an appropriate subsequence - that $x_n\asyleins$.
Let $(x_{n_\gamma})$ be a universal net. By Lemma \ref{lem_HP}
 it \wst-converges
to a point $x_s\in X_s$ and
\bglst
\lim_{\gamma}\norm{x+x_{n_\gamma}}\geq\norm{x+x_s}
=\norm{x}+\lim_{\gamma}\norm{x_{n_\gamma}}
\eglst
by \wst-continuity of the norm
which contradicts \Ref{glBem3} and proves property $(m_1^*)$.\\
(ii)$\gdw$(iii) is immediate from the definition of $\tau_{\mu}$.\\
(iii)$\folgt$(iv):
If $Z$ is separable then $\ball{X}$ is \wst-sequentially compact.
Hence $\ball{X}$ is $\tau_{\mu}$-sequentially compact if (iii) holds.\\
(iv)$\folgt$(i): If $Z$ is M-embedded then $Z=X^{\#}$. In this case, if
$\ball{X}$ is $\tau_{\mu}$-sequentially compact then the proof of Proposition
\ref{theo2.1b} shows that $Z$ has $(m_1^*)$.\\
(b) The Schur property clearly implies $(m_1)$.
Conversely, suppose a Banach space $X$ has $(m_1)$ but fails to have the
Schur property. Then there exists a normalized \w-null sequence $(x_n)$
in $X$.
Exactly as in (i)$\folgt$(ii) of part (a), using \Ref{glBem2} instead
of \Ref{glBem1}, we can extract a subsequence $(x_{n_k})$ which spans
$\leins$ asymptotically. But a normalized sequence spanning $\leins$ cannot
converge weakly. This contradiction proves that $(m_1)$ implies
the Schur property.
\ebew\bigskip\\
\begin{sublemma}\label{theo2.6}
For any L-embedded space $X$ one has the following inclusions
\bglst
X\cap \overline{\ball{X_s}}^{w^*}&\subset& \bigcap\{\overline{V}^w\mdE V\;\,\tau_{\mu}-\mbox{open in }
\ball{X},\,0\in V\}\\
&\subset& \bigcap\{\overline{\mbox{co}}^{\norm{\cdot}}(V)
\mdE V\;\,\tau_{\mu}-\mbox{open in } \ball{X},\,0\in V\}.
\eglst
\end{sublemma}
{\em Proof}:
The first inclusion is immediate from (i) of Lemma \ref{lem_elementar},
the second inclusion follows from
$\overline{\mbox{co}}^{w}(V)
=\overline{\mbox{co}}^{\norm{\cdot}}(V)$.\ebew\bigskip\\
That at the time of this writing we are not able to give the whole analogue
of \cite[L.\ 2.6]{GoKaLi} is essentially due to the fact that we know only
the behaviour of sequences but not of nets in $\tau_{\mu}$, more specifically if we knew
that one can extract a $\tau_{\mu}$-convergent sequence from a
$\tau_{\mu}$-convergent net then we could also prove that
$\bigcap\overline{\mbox{co}}^{\norm{\cdot}}(V)\subset X\cap \overline{\ball{X_s}}^{w^*}$.
In any case, our reduced version suffices for Lemma \ref{theo2.7}
which corresponds to \cite[L.\ 2.7]{GoKaLi}.
Note that (ii) of Lemma \ref{theo2.7} is equivalent to the two
conditions of Theorem \ref{theo3.10}.
\begin{lemma}\label{theo2.7}
Let $X$ be an L-embedded Banach space. The following assertions are equivalent.
\begin{description}
\item{\makebox[1.8em]{(i)}}
There exists a locally convex Hausdorff topology on $X$ which is coarser than
$\tau_{\mu}$ on $\ball{X}$.
\item{\makebox[1.8em]{(ii)}}
$\ball{X}$ is compact Hausdorff with respect to $\sigma(X,X^{\#})$.
\item{\makebox[1.8em]{(iii)}}
$\{0\}$ is the intersection of the convex $\tau_{\mu}$-neighbourhoods of $0$
in $\ball{X}$.
\end{description}
If moreover $\ball{X}$ is $\tau_{\mu}$-compact Hausdorff,
then the above conditions are also
equivalent to
\begin{description}
\item{\makebox[1.8em]{(iv)}}
The weak topology of $X$ is finer than $\tau_{\mu}$ on $\ball{X}$.
\item{\makebox[1.8em]{(v)}}
There exists a locally convex Hausdorff topology - namely $\sigma(X,X^{\#})$ -
which coincides with $\tau_{\mu}$ on $\ball{X}$.
\item{\makebox[1.8em]{(vi)}}
$0$ has a basis of $\tau_{\mu}$-neighbourhoods in $\ball{X}$ consisting
of convex sets.
\end{description}
\end{lemma}
{\em Proof}:
The proof is word-for-word the same as the one of \cite[L.\ 2.7]{GoKaLi}
(with Sublemma \ref{theo2.6} replacing \cite[L.\ 2.6]{GoKaLi}).
\ebew\bigskip\smallskip\\
Suppose that with the notation of (i) of Lemma \ref{lem_elementar}
one has $x=x''\in X$; then $x_{\gamma}\gentaumu0$ and $x_{\gamma}\genw x$.
If (iv) of Lemma \ref{theo2.7} holds then this means that both
$x_{\gamma}\gentaumu 0$ and $x_{\gamma}\gentaumu x$.
Since it is not clear whether $\tau_{\mu}$ is
Hausdorff  we cannot deduce that $x=0$. Therefore we claimed $\ball{X}$ to
be $\tau_{\mu}$-Hausdorff for (iv) - (vi) in Lemma \ref{theo2.7}
(in particular for the proof of (iv)$\folgt$(ii)).
\bigskip\\
We end this section with some remarks on the continuity of the canonical
L-projection $P$ on the bidual of an L-embedded space $X$.

In \cite[Prop.\ 3.8]{GoKaLi} it was shown that $P$ is
($w^*$-$\tau_{\mu}$)-sequentially continuous if and only if $X$ has the Schur property.
In the proof the authors of \cite{GoKaLi} use the Grothendieck property of
$\Lunendl$ in  order to obtain the ($w^*$-$w$)-sequential continuity of the
L-projection  on $(\Lunendl)'$.
But in the general case it is not known whether $P$ is always
($w^*$-$w$)-sequentially continuous although by an example of Johnson
it is known that in general the dual of an L-embedded Banach space does not
have the Grothendieck property.
Therefore it is not clear whether the first part of
\cite[Prop.\ 3.8]{GoKaLi} can be
generalized in our context; its second part can:
\begin{proposition}\label{theo3.8b}
Let $X$ be L-embedded and $P$ be the canonical L-projection on $X''$.
Suppose that the abstract measure topology $\tau_{\mu}$ on $X$ is
Hausdorff.
Then $P$ is ($w^*$-$\tau_{\mu}$)-continuous if and only if the restriction of
$\tau_{\mu}$ to $\ball{X}$ is compact Hausdorff locally convex.
\end{proposition}
{\em Proof}: See \cite[Prop.\ 3.8(b)]{GoKaLi}.\ebew\bigskip\bigskip\\
{\bf\S 6 The special case of the predual of a von Neumann algebra}\\
Recall that two elements $\phi, \psi\in{\cal N}_*$ in the predual of
a von Neumann algebra ${\cal N}$ are called orthogonal - $\phi\bot\psi$
in symbols - if
they have orthogonal right and orthogonal left support projections.
It is well know that $\phi\bot\psi$ if and only if the linear span
of $\phi$ and $\psi$ is isometrically isomorphic to the two-dimensional
$\leins(2)$.
\begin{sublemma}\label{sublem_vN}
Let ${\cal N}$ be a von Neumann algebra with predual $X={\cal N}_*$.
For every $\eps>0$ there is $\delta>0$ with the following property.\\
If $x,y,z\in \ball{X}$ are such that 
\bglst
\Norm{\alpha\frac{z}{\norm{z}}+\beta\frac{x}{\norm{x}}}
&\geq&
(1-\delta)\Bigl(\betr{\alpha}+\betr{\beta}\Bigr)\\
\Norm{\alpha\frac{z}{\norm{z}}+\beta\frac{y}{\norm{y}}}
&\geq&
(1-\delta)\Bigl(\betr{\alpha}+\betr{\beta}\Bigr)
\eglst
then there are $\tilde{x}, \tilde{y}, \tilde{z}\in X$ of norm one
such that
\bglst
\begin{array}{r}
\tilde{x}\bot\tilde{z}\\
\tilde{y}\bot\tilde{z}
\end{array}
\mbox{ and }
\begin{array}{r}
\norm{x-\tilde{x}}\leq\eps\\
\norm{y-\tilde{y}}\leq\eps\\
\norm{z-\tilde{z}}\leq\eps\\
\end{array}
\eglst
\end{sublemma}
{\em Proof}:
We have already recalled the fact that two normal functionals on a
von Neumann algebra are orthogonal if and only if they span $\leins(2)$
isometrically. A second ingredient of this proof is the fact 
that the ultraproduct of a family of preduals of von Neumann
algebras is again a predual of a von Neumann algebra
\cite{Groh}, \cite{preRaynaud}.
It is now enough to combine these two facts with
a standard ultraproduct argument \cite{hei1}.
\ebew\bigskip\\It is well known that if ${\cal N}$ is a von Neumann algebra with a
finite faithful normal trace $\tau$ then
the measure topology on the predual ${\cal N}_*$ defined via $\tau$ is metrizable
and makes ${\cal N}_*$ a Hausdorff topological vector space.
For the abstract measure topology $\tau_{\mu}$ on the predual of an arbitrary
von Neumann algebra we only know by the following lemma (and by Lemmas
\ref{lem_elementar}, \ref{lem_Frechet}) that
multiplication by scalars and addition are continuous and that
on bounded sets $\tau_{\mu}$ makes the von Neumann predual
a Fr\'echet space.
\begin{proposition}\label{theo_add}
Let $X={\cal N}_*$ be the predual of a von Neumann algebra ${\cal N}$.
If $X$ is endowed with the abstract measure toplogy $\tau_{\mu}$
which it has by Theorem \ref{theo_existence} then addition is
$\tau_{\mu}$-continuous.
\end{proposition}
{\em Proof}:
Let $(x_n)$, $(y_n)$ be two sequences in $\ball{X}$ each of which
spans $\leins$ asymptotically.
Suppose there is $\eps>0$ such that $\norm{z_n}>\eps$ where $z_n=x_n+y_n$,
suppose further that $\lim\norm{x_n}$, $\lim\norm{y_n}$, $\lim\norm{z_n}$
exist.
To show the Proposition it is enough to show that there is a subsequence
$(z_{n_k})$ spanning $\leins$ asymptotically.

Let $(\delta_n)$ be a sequence of strictly positive numbers
in $]0,1]$ converging to $0$.
Set $\eta_1=\frac{1}{3}\delta_1$ and
$\eta_{n+1}=\frac{1}{3}\min(\eta_n,\delta_{n+1})$ for
$n\in\N$.
By induction over $m\in\N$ we construct a sequence $(n_m)$
such that 
\bgl
\Norm{\sum_{k=1}^{m}\alpha_k\frac{z_{n_k}}{\norm{z_{n_k}}}}
\geq
\Bigl(\sum_{k=1}^{m}(1-\delta_k)\betr{\alpha_k}\Bigr)
+\eta_m \sum_{k=1}^{m}\betr{\alpha_k}
\label{glvN1}
\egl
for all integers $m$ and all scalars $\alpha_k$.
Since the construction ressembles the proofs of \cite[Th.\ 2]{Pfi-Fixp},
of Lemma \ref{lem_Frechet}, and of Remark \ref{Bem} we detail it only
for $m=2$.

Let $\theta_2>0$ be arbitrary for the moment and to be determined later.
Set $n_1=1$.\\
{\sc Claim}:
There is $n_2\in\N$ such that 
\bgl
\Norm{\alpha\frac{z_1}{\norm{z_1}}+\beta\frac{x_{n_2}}{\norm{x_{n_2}}}}
&\geq&
(1-\theta_2)\Bigl(\betr{\alpha}+\betr{\beta}\Bigr)\label{gladd1}\\
\Norm{\alpha\frac{z_1}{\norm{z_1}}+\beta\frac{y_{n_2}}{\norm{y_{n_2}}}}
&\geq&
(1-\theta_2)\Bigl(\betr{\alpha}+\betr{\beta}\Bigr)
\egl
for all scalars $\alpha, \beta$.
(As usual in this paper) this is essentially due to the
\wst-lower semicontinuity of the norm by which one has
$\liminf_{\gamma}
\Norm{\alpha\frac{z_1}{\norm{z_1}}+\beta\frac{x_{n_\gamma}}{\norm{x_{n_\gamma}}}}
\geq
\Norm{\alpha\frac{z_1}{\norm{z_1}}+\beta\frac{x_s}{\norm{x_s}}}
=\betr{\alpha}+\betr{\beta}$
where $x_s\in X_s$ is the \wst-limit of a universal net $(x_{n_\gamma})$
and $\norm{x_s}=\lim_{\gamma}\norm{x_{n_\gamma}}$ (see Lemma \ref{lem_HP});
hence there are infinitely many $n_2$ satisfying \Ref{gladd1}.
Applying the same reasoning to the corresponding subsequence of
$(y_n)$ gives the {\sc Claim}.

By Sublemma \ref{sublem_vN} there are
$\tilde{x}_2$, $\tilde{y}_2$, $\tilde{z}_1$ in $X$ such that
\bglst
\norm{x_{n_2}-\tilde{x}_2}\leq \theta'_2, \mbox{ }
\norm{y_{n_2}-\tilde{y}_2}\leq \theta'_2, \mbox{ }
\norm{z_1-\tilde{z}_1}\leq \theta'_2,
\eglst
and
\bgl
\tilde{z}_1 \,\bot\, \tilde{x}_2, \mbox{ }\mbox{ }
\tilde{z}_1 \,\bot\, \tilde{y}_2   \label{glvN2}
\egl
where $\theta'_2 \gen0$ as $\theta_2 \gen0$.
Now \Ref{glvN2} implies that
\bglst
\tilde{z}_1 \,\bot\, \tilde{z}_{n_2}
\eglst
where $\tilde{z}_{n_2}=\tilde{x}_{n_2} +\tilde{y}_{n_2}$
which means that
$\Norm{\alpha\frac{\tilde{z}_1}{\norm{\tilde{z}_1}}
+\beta\frac{\tilde{z}_2}{\norm{\tilde{z}_2}}}
=\betr{\alpha}+\betr{\beta}$.
Hence, if $\theta_2$ is choosen small enough, we get \Ref{glvN1} for $m=2$.
It is now left to the reader to iterate the construction in order to end the induction and
thus the proof.
\ebew\bigskip\medskip\\
As already noticed in the remarks concerning the questions after Corollary
\ref{coro_compakt}, Randrianantoanina \cite{Randri-Kad-Pel}
has proved that the Kadec-\Pel\ splitting lemma
holds for preduals of von Neumann algebras.
With this result Komlos' theorem follows almost immediately
for von Neumann preduals.
\begin{proposition}[Komlos]\label{theoKomlosvN}
Each bounded sequence in the predual of a von Neumann algebra admits a
subsequence whose Cesaro means converge with respect to the
abstract measure topology.
\end{proposition}
{\em Proof}:
Let $X={\cal N}_*$ be the predual of a von Neumann algebra ${\cal N}$.
Endow $X$ with its abstract measure topology $\tau_{\mu}$.
Let $(x_n)\subset X$ be bounded.
Then by \cite{Randri-Kad-Pel} there is a subsequence $(x_{n_k})$ and there
is a decomposition $x_{n_k}=y_k + z_k$ where
$(z_k)$ \w-converges and such that there is a sequence $(\tilde{y}_k)$
of pairwise orthogonal elements in $X$ with
\bgl
\norm{\tilde{y}_k -y_k}\gen 0.\label{glKomvN1}
\egl
If one applies the classical theorem of Komlos to the
isometric $\leins$-copy  $[\tilde{y}_k]$
then one obtains, after passing to an appropriate
subsequence, that the Cesaro means of
any subsequence of $(\tilde{y}_k)$ converge with respect to the measure
(=pointwise) topology of $[\tilde{y}_k]$ whence with respect to $\tau_{\mu}$.
So do the Cesaro means $\hat{y}_k$ of the $y_k$ by \Ref{glKomvN1}.
It is known that von Neumann preduals have the weak Banach-Saks property
\cite{BelDie}.
Hence, again after passing to the appropriate subsequences of $(y_k)$ and
$(z_k)$,
the Cesaro means $\hat{z}_k$ of the $z_k$ converge in norm
whence with respect to $\tau_{\mu}$.
Since $\hat{x}_{k}=\hat{y}_k +\hat{z}_k$ where $\hat{x}_{k}$ denote
the Cesaro means of the $x_{n_k}$, the assertion now follows from Proposition
\ref{theo_add}.
\ebew\bigskip\medskip\\
{\sc Acknowledgement}
I thank Dirk Werner for several helpful discussions.

\begin{thebibliography}{10}

\bibitem{Beauz}
B.~Beauzamy.
\newblock {Banach-Saks properties and spreading models}.
\newblock {\em Math. Scand.}, 44:357--384, 1979.

\bibitem{BelDie}
A.~Belanger and J.~Diestel.
\newblock {A remark on weak convergence in the dual of a
  C{$^{*}$}-al\-ge\-bra}.
\newblock {\em Proc. Amer. Math. Soc.}, 98:185--186, 1986.

\bibitem{Buh-Lo}
A.~V. Buhvalov and G.~J. Lozanovskii.
\newblock On sets closed in measure in spaces of measurable functions.
\newblock {\em Trans. Moscow Math. Soc.}, 2:127--148, 1978.

\bibitem{Die-Seq}
J.~Diestel.
\newblock {\em Sequences and Series in Banach Spaces}.
\newblock Springer, Berlin-Heidelberg-New York, 1984.

\bibitem{Dor}
L.~E. Dor.
\newblock {On projections in {$L_1$}}.
\newblock {\em Ann. of Math.}, 102:463--474, 1975.

\bibitem{DJLT}
P.~N. Dowling, W.~B. Johnson, C.~J. Lennard, and B.~Turett.
\newblock {The optimality of James's distortion theorems}.
\newblock {\em Proc. Amer. Math. Soc.}, 125:167--174, 1997.

\bibitem{DS}
N.~Dunford and J.~T. Schwartz.
\newblock {\em Linear Operators. Part 1: General Theory}.
\newblock Interscience Publishers, New York, 1958.

\bibitem{Engel2}
R.~Engelking.
\newblock {\em General Topology}.
\newblock {Heldermann} {Verlag}, {Berlin}, 1989.

\bibitem{G-bien}
G.~Godefroy.
\newblock Sous-espaces bien dispos\'{e}s de {$L^{1}$} -- {Applications}.
\newblock {\em Trans. Amer. Math. Soc.}, 286:227--249, 1984.

\bibitem{GoKaLi}
G.~Godefroy, N.~Kalton, and D.~Li.
\newblock On subspaces of {$L^1$} which embed into {$l_1$}.
\newblock {\em J. Reine Angew. Math.}, 471:43--75, 1996.

\bibitem{G-Li}
G.~Godefroy and D.~Li.
\newblock Some natural families of {$M$}-ideals.
\newblock {\em Math. Scand.}, 66:249--263, 1990.

\bibitem{Groh}
U.~Groh.
\newblock Uniform ergodic theorems for identity preserving {Schwartz} maps on
  {$W^*$}-algebras.
\newblock {\em J. Operator Theory}, 11:395--404, 1984.

\bibitem{HWW}
P.~Harmand, D.~Werner, and W.~Werner.
\newblock {\em {$M$}-ideals in {Banach} {Spaces} and {Banach} {Algebras}}.
\newblock Lecture Notes in Mathematics 1547. Springer, 1993.

\bibitem{hei1}
S.~Heinrich.
\newblock Ultraproducts in {Banach} space theory.
\newblock {\em J. Reine Angew. Math.}, 313:72--104, 1980.

\bibitem{KadPel}
I.~Kadec and A.~Pe{\l}czy\'{n}ski.
\newblock {Bases, lacunary sequences and complemented subspaces in the space
  ${\mbox{L}}_{p}$.}
\newblock {\em Studia Math.}, 21:161--176, 1962.

\bibitem{KaWe}
N.~Kalton and D.~Werner.
\newblock Property ({$M$}), {$M$}-ideals, and almost isometric structure of
  {Banach} spaces.
\newblock {\em J. Reine Angew. Math.}, 461:137--178, 1995.

\bibitem{Kelley}
J.~L. Kelley.
\newblock {\em {{General} {Topology}}}.
\newblock {D. van Nostrand Company}, 1963.

\bibitem{Kisynski}
J.~Kisy\'{n}ski.
\newblock {{Convergence} du type {L}}.
\newblock {\em Colloquium Mathematicum}, 7:205--211, 1960.

\bibitem{Li-Ox}
D.~Li.
\newblock Espaces {$L$}-facteurs de leurs biduaux: bonne disposition, meilleure
  approximation et propri\'{e}t\'{e} de {Radon-Nikodym}.
\newblock {\em Quart. J. Math. Oxford (2)}, 38:229--243, 1987.

\bibitem{LiTz1}
J.~Lindenstrauss and L.~Tzafriri.
\newblock {\em Classical Banach Spaces I}.
\newblock Springer, Berlin-Heidelberg-New York, 1977.

\bibitem{LiTz2}
J.~Lindenstrauss and L.~Tzafriri.
\newblock {\em Classical Banach Spaces II}.
\newblock Springer, Berlin-Heidelberg-New York, 1979.

\bibitem{Ped}
G.~K. Pedersen.
\newblock {\em {C{$^{*}$}-algebras and Their Automorphism Groups}}.
\newblock Academic Press, London, New York, San Francisco, 1979.

\bibitem{Pfi-Fixp}
H.~Pfitzner.
\newblock {A note on asymptotically isometric copies of {$l^1$} and {$c_0$}}.
\newblock {\em Proc. Amer. Math. Soc.}, 2000.
\newblock to appear.

\bibitem{Pfi-vN}
H.~Pfitzner.
\newblock {Perturbations of {$l^1$}-copies and convergence in preduals of von
  {Neumann} algebras}.
\newblock {\em J. Operator Theory}, 2000.
\newblock to appear.

\bibitem{Pi-bases}
G.~Pisier.
\newblock {Bases, suites lacunaires dans les espaces ${\mbox{L}}_{p}$
  d'apr\`{e}s Kadec et Pelczynski.}
\newblock {\em {S}\'{e}minaire Maurey-Schwartz de l'{E}cole {P}olytechnique},
  pages Expos\'{e}s 18 -- 19, 1972-73.

\bibitem{Randri-Kad-Pel}
N.~Randrianantoanina.
\newblock {Kadec-Pe\l czy\'nski} decomposition for {Haagerup} {$L^{p}$}-spaces.
\newblock Preprint, 2000.

\bibitem{preRaynaud}
Yves Raynaud.
\newblock On ultrapowers of non commutative {$L_p$} spaces.
\newblock Preprint, 2000.

\bibitem{Tak}
M.~Takesaki.
\newblock {\em Theory of Operator Algebras {I}}.
\newblock Springer, Berlin-Heidelberg-New York, 1979.

\end{thebibliography}

Hermann Pfitzner\\
Universit\'e d'Orl\'eans\\
BP 6759\\
F-45067 Orl\'eans Cedex 2\\
France\\
e-mail: pfitzner@labomath.univ-orleans.fr
\end{document}